\documentclass[11pt,a4paper]{amsart}

\usepackage[T1]{fontenc}
\usepackage[a4paper,margin=25mm]{geometry}
\usepackage{amsmath,amssymb,mathtools}
\usepackage{enumitem}
\usepackage{microtype}

\usepackage[hidelinks]{hyperref}
\hypersetup{
  pdftitle={Kac lemma on noncommutative probability spaces},
  pdfauthor={Tattwamasi Amrutam, Diptesh Saha, and Xumin Wang},
  pdfsubject={Kac-type formulas for actions on finite von Neumann algebras},
  pdfkeywords={Kac lemma, von Neumann algebra, countable group action,
    partial isometry, return-time observable, center-valued trace},
  pdfdisplaydoctitle=true
}
\allowdisplaybreaks
\numberwithin{equation}{section}
\setlist[enumerate]{label=\textup{(\roman*)},leftmargin=2.4em}

\newtheorem{theorem}{Theorem}[section]
\newtheorem{proposition}[theorem]{Proposition}
\newtheorem{lemma}[theorem]{Lemma}
\newtheorem{corollary}[theorem]{Corollary}
\theoremstyle{definition}
\newtheorem{definition}[theorem]{Definition}
\newtheorem{example}[theorem]{Example}
\theoremstyle{remark}
\newtheorem{remark}[theorem]{Remark}

\newcommand{\M}{\mathcal M}
\newcommand{\C}{\mathbb C}
\newcommand{\Z}{\mathbb Z}
\newcommand{\N}{\mathbb N}
\newcommand{\one}{\mathbf 1}
\newcommand{\Proj}{\operatorname{Proj}}
\newcommand{\Aut}{\operatorname{Aut}}
\newcommand{\Tr}{\operatorname{Tr}}
\newcommand{\abs}[1]{\left\vert#1\right\vert}

\title[Kac lemma on noncommutative probability spaces]
  {Kac lemma on noncommutative probability spaces}
\author[Amrutam]{Tattwamasi Amrutam}
\address{Institute of Mathematics of the Polish Academy of Sciences, ul. Sniadeckich 8, 00-656, Warszawa, Poland}
\email{tattwamasiamrutam@gmail.com}
\thanks{T.A. is supported by National Science Centre, Poland Sonata, grant number 2025/59/D/ST1/03117.}

\author[Saha]{Diptesh Saha}
\address{Institute of Mathematics of the Polish Academy of Sciences, ul. Sniadeckich 8, 00-656, Warszawa, Poland\\\text{and,}\newline
Indian Statistical Institute, 203, Barrackpore Trunk Road
Kolkata, West Bengal 700108, India}
\email{dptshs@gmail.com}
\author[Wang]{Xumin Wang}
\address{Institute for Advanced Study in Mathematics, Harbin Institute of
Technology, Harbin, 150001, China}
\thanks{X.W. is supported by the National Natural Science Foundation of China, grant numbers 12601242 and W2611005. }
\email{xumin.wang@hit.edu.cn}
\date{}

\begin{document}

\begin{abstract}
We prove non-commutative Kac's formula for actions of countable groups. We replace the return time functions by suitable projections in von Neumann algebras. Moreover, we prove a generalised Kac's formula in the non-commutative setting which extends the allocation construction of \cite{mw24} and is motivated by the random Kac inequality of \cite{aaronson1996az}.
\end{abstract}

\maketitle

\section{Introduction}

Poincare's celebrated recurrence theorem \cite{brown1976ergodic} states that for a probability space $(X, \mathcal F, \mu)$, a probability preserving transformation $T$ on $X$, and a set $E \in \mathcal F$ with $\mu(E)>0$, for almost every $x \in E$, we have $T^n(x) \in E$ for some $n \in \mathbb N$. Kac's lemma \cite{kac1947} studies the expected time taken by a point to return to its set of origin in a probability space with an ergodic probability preserving transformation defined on it. Moreover, it calculates the expected time to be the inverse of the probability measure of the set. The classical Kac lemma states the following.

For a probability preserving transformation $T$ and a measurable set $A$ of positive measure, let $r_A(x):=\min\{n\in\mathbb N:T^n(x)\in A\}$ denote the return time of $x\in A$, with $\min\varnothing:=\infty$.

\begin{theorem}\label{kac lem}
	Let $(X, \mathcal F, \mu)$ be a probability space and $T:X \to X$ be an ergodic measure preserving transformation. Let $A \in \mathcal F$ such that $\mu(A)>0$. Then $\int_A r_A(x)d\mu(x)=1$.
\end{theorem}

An extension of Theorem~\ref{kac lem} for ergodic probability preserving actions of the group $\mathbb Z^d$ was first studied by Aaronson and Weiss \cite{aaronson1996az}. For higher-rank actions, such as \(\mathbb Z^d\)-actions, there is no canonical notion of first return. Using the geometry of lattice boxes, they proved a random Kac inequality \cite[Theorem~2.1]{aaronson1996az}.

In \cite{mw24}, the authors replace the Kac function of Aaronson and Weiss \cite{aaronson1996az} by an allocation construction, and extend the result to arbitrary countable groups. In particular, for $\mathbb Z^d$, their result strengthens the integrability conclusion of Aaronson–Weiss, although it does not preserve the specific dyadic geometry of the \(\mathbb Z^d\)-boxes used in the Aaronson–Weiss formulation.  More precisely, their main result \cite[Theorem 1.5]{mw24} states that for every countable group $\Gamma$, there is a sequence of finite subsets $(B_n)_n$ of $\Gamma$ containing the identity $e$ such that for every probability preserving action and every sweep-out set $A$, there exists a measurable function $\phi:A \to \mathbb N$ so that the translates associated to $\phi$ give a measurable partition of $X$ modulo null sets, and $\int_A\abs{B_{\phi(x)}}~d\mu(x)=1$.
This identity is a conservation of mass formulation of Kac's lemma.

One of the earliest links between ergodic theory and operator algebras is provided by von Neumann's work on ergodic theorems and his work on the structure of operator algebras. Typically, a classical probability measure space $(X, \mathcal B, \mu)$ is replaced by a von Neumann algebra $M$ in the non-commutative setting, where a fixed state $\varphi$ plays the role of the probability measure. Furthermore, a state preserving automorphism on the von Neumann algebra is often considered as the non-commutative analogue of a probability measure preserving transformation and measurable sets are often replaced by projections.

From this perspective, Kac's lemma asks for a quantitative relation between the size of a region and the time required for the dynamics to return to it. In a noncommutative probability space there are no points and, consequently, no pointwise return-time functions in the classical sense. The correct replacement must therefore be expressed in terms of projections, states (or traces), and the action of the automorphism. The challenge is to identify an operator-algebraic object that plays the role of the classical return blocks or allocation cells and retains the fundamental conservation principle behind Kac's formula.

Non-commutative recurrence theorems were studied in \cite{nsz03}, while \L{}uczak \cite{luczak2005} studied time projections for stopping times in von Neumann algebras. Kac-type formulas for iterated open quantum systems were obtained in \cite{sinkovicz2016}. In this article, we develop a non-commutative analogue of Kac's formula for actions on finite von Neumann algebras. Our result is motivated by the classical results of Kac \cite{kac1947}, Aaronson, Weiss, Meyerovitch (\cite{aaronson1996az}, \cite{mw24}) and the non-commutative stopping-time framework of \L{}uczak \cite{luczak2005}.

We begin with a finite von Neumann algebra $\M$ with a faithful, normal state $\varphi$ and a group of $*$-automorphisms on $\M$ preserving $\varphi$. First we observe that if the group $G$ further acts on $\M$ ergodically, then $\varphi$ is a tracial state. So we will always consider a non-commutative dynamical system $(\M, \varphi)$, where $\M$ is a finite von Neumann algebra and $\varphi$ is a faithful normal tracial state. We refer to Section~\ref{sec:preliminaries} for relevant definitions and results. In Section~\ref{sec:Z-action}, we consider an ergodic $\varphi$-preserving action of $\mathbb Z$ and a non-zero projection $p\in\M$. We identify the non-commutative analogue of the event of a point returning to a set by the appropriate sequence of projections $(f_n)_n$ and define a non-commutative return time associated to a projection $p$ as a positive, self-adjoint operator $T_p:=\sum_{n\in\N}nf_n$ and show that $\varphi(T_p)=1$ and $\varphi_p(T_p)=1/\varphi(p)$, where $\varphi_p=\varphi(p)^{-1}\varphi|_{p\M p}$ (see Theorem \ref{thm:quantum-return-time}). The proof uses projection lattice techniques to handle non-commutativity.

In Section~\ref{sec:universal-towers}, we generalise Meyerovitch and Weiss's allocation construction \cite{mw24} for general countable group actions to the non-commutative setting. Due to the noncommutative nature of the underlying von Neumann algebra, the projections corresponding to the classical events of a point returning to a set and its image under an automorphism need not satisfy an exact inclusion or equality. The natural noncommutative analogue is therefore expressed in terms of Murray–von Neumann subequivalence. This can obstruct the existence of a direct $p$-allocation (the non-commutative analogue of a classical allocation function). As a result, we work with weak partial-isometric allocations (see Definition \ref{def:partial-isometric-allocation}) and provide an example (see Example \ref{ex:direct-allocation-fails}) where a direct $p$-allocation does not exist. Finally, we prove the following non-commutative Kac formula for countable groups generalising Theorem \ref{thm:quantum-return-time}.

For a partial-isometric $p$-allocation $(\rho_\gamma)_{\gamma\in\Gamma}$ with transport witness $(p_\gamma,u_\gamma)_{\gamma\in\Gamma}$ and $x\in\M$, consider the series $L_{\rho,p,u}(x):=\sum_{\gamma\in\Gamma}u_\gamma^*\alpha_\gamma(x)u_\gamma$.

\begin{theorem}
\label{thm:order-free-allocation-intro}
Let $\alpha\colon\Gamma\to\Aut(\M)$ be an action of a countable group preserving the faithful normal tracial state $\varphi$ on the finite von Neumann algebra $\M$, and let $p\in\M$ be a non-zero projection.
Let \((\rho_\gamma)_{\gamma \in \Gamma}\) be a partial-isometric \(p\)-allocation with transport witness \((p_\gamma,u_\gamma)_{\gamma \in \Gamma}\).  For every \(x\in\M\), the series $L_{\rho,p,u}(x)$
converges absolutely in
\( L^1(p\M p,\varphi)\). The sum satisfies
\begin{enumerate}
    \item \(\sum_{\gamma\in\Gamma}
\|u_\gamma^*\alpha_\gamma(x)u_\gamma\|_1\leq\|x\|\) and
\item \(\varphi(  L_{\rho,p,u}(x))=\varphi(x)\).
\end{enumerate}
Moreover, if \(x\geq0\), this limit is also the supremum of the increasing net
of finite subsums in the extended positive cone of \(p\M p\).
\end{theorem}

In \cite[Theorem 1.5]{mw24}, the authors prove the existence of a universal sequence of finite subsets $(B_n)_n$ of the group $\Gamma$ and, for every sweep-out set $A$, construct a measurable function $\phi:A\to\mathbb N$ such that the corresponding translates of the level sets $\phi^{-1}(n)$ partition the set $X$. Generalising this result, we construct universal partial-isometric Kac towers. Given a finite von Neumann algebra $(\M,\tau)$ with a faithful normal tracial state, a $\tau$-preserving action $\alpha\colon\Gamma\to\Aut(\M)$ and a non-zero sweep-out projection $p\in\M$, Proposition~\ref{prop:universal-kac-towers1} gives a sequence $(B_n)_n$ of finite subsets of $\Gamma$ depending only on $\Gamma$, projections $(q_n)_n$, $(\rho_{n,\gamma})_{n,\gamma}$ and partial isometries $(v_{n,\gamma})_{n,\gamma}$ which form a partial-isometric Kac tower with shapes $(B_n)_n$.

To construct the partial-isometric Kac tower in Proposition~\ref{prop:universal-kac-towers1}, we first fix a sequence of finite subsets $(B_n)$ of $\Gamma$ which contains every finite subset of $\Gamma$ containing the identity. For the given trace-preserving action and sweep-out projection $p$, Lemma~\ref{lem:greedy-quantum-allocation} gives source projections and subprojections of $p$. Applying Lemma~\ref{lem:finite-shape-threshold-coupling} to the center-valued traces of these subprojections yields a family of central positive elements indexed by the finite sets, whose sums recover those traces.

These central elements satisfy the type~I integrality condition of Corollary~\ref{cor:center-valued-dimension-partition}, which allows them to be realized by mutually orthogonal projections. After refining the source projections according to the elements of each $B_n$, equivariance of the center-valued trace shows that the corresponding target projections have matching center-valued dimensions. The center-valued dimension criterion therefore produces partial isometries identifying each base projection with the corresponding translate of a source projection. Finally, reindexing by the fixed sequence $(B_n)$ yields the required family.

Let $(\M,\tau)$ be a finite von Neumann algebra with a faithful normal tracial state, and let $\alpha\colon\Gamma\to\Aut(\M)$ be a $\tau$-preserving action of a countable discrete group $\Gamma$. Let $p\in\M$ be a non-zero projection satisfying $\bigvee_{\gamma\in\Gamma}\alpha_\gamma(p)=1$. Fix a partial-isometric Kac tower $(q_n,\rho_{n,\gamma},v_{n,\gamma})$ for $(\alpha,p)$ with shapes $(B_n)_{n\geq1}$, in the sense of Definition~\ref{def:kac-tower}.

For $x\in L^1(\M,\tau)$, consider the series $L_p(x):=\sum_{n\geq1}\sum_{\gamma\in B_n}v_{n,\gamma}^*\alpha_\gamma(x)v_{n,\gamma}$. Let \(R_p:=\sum_{n\geq1}|B_n|q_n\) and \(\Phi_p:=\sum_{n\geq1}nq_n\), where both sums are positive spectral series. For \(n\geq1\), let \(h(n):=|B_n|\), and let \(\tau_p(x):=\tau(x)/\tau(p)\) be the normalized corner trace. We prove the following generalised Kac formula.

\begin{theorem}\label{thm:main thm-intro}
For every $x\in L^1(\M,\tau)$, the series $L_p(x)$ converges absolutely in
\( L^1(p\M p,\tau|_{p\M p})\). Moreover,
we have $\sum_{n\geq1}\sum_{\gamma\in B_n}\|v_{n,\gamma}^*\alpha_\gamma(x)v_{n,\gamma}\|_1\leq\|x\|_1$ and $\tau(L_p(x))=\tau(x)$.
Furthermore, \(R_p=L_p(1)\) belongs to
\( L^1(p\M p,\tau|_{p\M p})\) and satisfies
\(\tau(  R_p)=\sum_{n\geq1}|B_n|\tau(q_n)=1\).
Consequently, \(\tau_p(R_p)=1/\tau(p)\).

Moreover, the positive operator \(\Phi_p\) is affiliated
with \(p\M p\). For each \(n\geq1\), we have \(\one_{\{n\}}(\Phi_p)=q_n\). We also have
\(R_p=h(\Phi_p)\).
\end{theorem}
\subsection*{Acknowledgements}
 This work was partially supported by the Simons Foundation grant (award no. SFI-MPS-T-Institutes-00010825) and from State Treasury funds as part of a task commissioned by the Minister of Science and Higher Education under the project “Organization of the Simons Semesters at the Banach Center - New Energies in 2026-2028” (agreement no. MNiSW/2025/DAP/491). 
\vspace{0.5cm} % Adds a little breathing room before the logo
\begin{center}
    % Adjust the 0.8 to make the logo wider or narrower relative to the text width
    \includegraphics[width=0.5\textwidth]{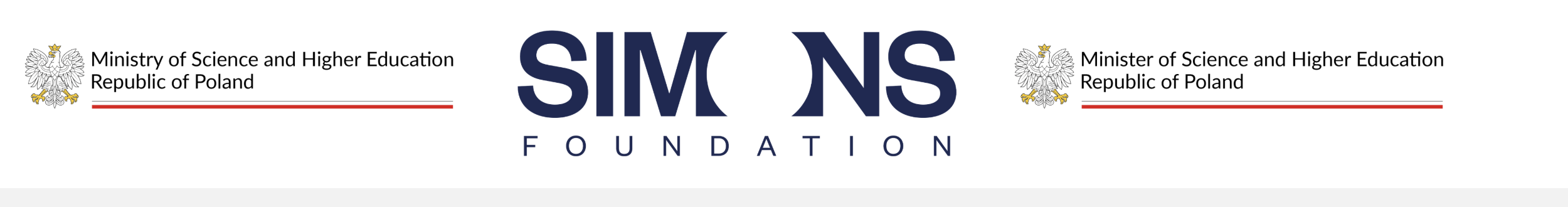} 
\end{center}
\section{Preliminaries}
\label{sec:preliminaries}

Let us fix a finite von Neumann algebra \(\M\) equipped with a faithful
normal state \(\varphi\).  For an action \(\alpha:G\to\Aut(\M)\) and
\(g,h\in G\), we use the convention \(\alpha_{gh}=\alpha_g\circ\alpha_h\). Thus \(\alpha_g^{-1}=\alpha_{g^{-1}}\).  We say that
\(\alpha\) preserves \(\varphi\) if \(\varphi\circ\alpha_g=\varphi\)
for every \(g\in G\).  The fixed-point algebra \(\M^\alpha\) consists
of those \(x\in\M\) for which \(\alpha_g(x)=x\) for every \(g\in G\),
and the action is called ergodic if \(\M^\alpha=\C1\).  We denote by
\(\M_\varphi\) the centralizer of \(\varphi\); thus \(a\in\M_\varphi\)
if and only if \(\varphi(ay)=\varphi(ya)\) for every \(y\in\M\).

We first observe that an invariant faithful normal state on a finite
von Neumann algebra is tracial whenever the action is ergodic.

\begin{proposition}
\label{prop:invariant-state-is-trace}
Let \(G\) act ergodically and \(\varphi\)-preservingly on the finite
von Neumann algebra \(\M\).  Then \(\varphi\) is a faithful normal
tracial state.  In particular, \(\M_\varphi=\M\).
\end{proposition}

\begin{proof}
Let \(\Theta:\M\to Z(\M)\) be the normalized center-valued trace
\cite[Theorem~III.2.5.7]{MR2188261}, and let
\(\tau=\varphi\circ\Theta\).  Since \(\Theta\) is faithful, normal,
unital and tracial, \(\tau\) is a faithful normal tracial state.
For each \(g\in G\), the map
\(\alpha_{g^{-1}}\circ\Theta\circ\alpha_g\) is again a normalized
center-valued trace on \(\M\).  By uniqueness, this map equals
\(\Theta\), and hence
\(\Theta\circ\alpha_g=\alpha_g\circ\Theta\).  It follows that
\(\tau\circ\alpha_g=\varphi\circ\alpha_g\circ\Theta
=\varphi\circ\Theta=\tau\).

By the identification of the positive predual with
\( L^1(\M,\tau)_+\), there is a unique positive element
\(h\in  L^1(\M,\tau)\) such that, for every \(x\in\M\), we have
\(\varphi(x)=\tau(hx)\) \cite[Corollary~4.51]{hiai2021}.  Each
\(\tau\)-preserving automorphism \(\alpha_g\) extends to
\( L^1(\M,\tau)\).  For \(x\in\M\) and \(g\in G\), invariance of
\(\tau\) and \(\varphi\) gives
\[
  \tau\bigl(\alpha_{g^{-1}}(h)x\bigr)
  =\tau\bigl(h\alpha_g(x)\bigr)
  =\varphi(\alpha_g(x))
  =\varphi(x)
  =\tau(hx).
\]
Thus \(\alpha_{g^{-1}}(h)\) is another positive density for
\(\varphi\), and uniqueness gives \(\alpha_{g^{-1}}(h)=h\).  Since
\(g\) was arbitrary, \(h\) is fixed by every \(\alpha_g\).
Consequently, every spectral projection of \(h\) belongs to
\(\M^\alpha=\C1\), and so each such projection is either \(0\) or
\(1\).  The spectral resolution of the positive integrable operator
\(h\) therefore gives \(h=c1\) for some \(c\geq0\).  Finally,
\(1=\varphi(1)=\tau(h)=c\), so \(\varphi=\tau\).
\end{proof}

Under these hypotheses we retain the notation \(\varphi\) for the
trace and let \(\|x\|_1=\varphi(|x|)\).  We use
\(L_0(\M,\varphi)\), \( L^1(\M,\varphi)\), and the measure topology
on \(L_0(\M,\varphi)\) in the sense of
\cite[Sections~4.1--4.3]{hiai2021}.  Since \(\varphi(1)<\infty\),
every closed densely defined operator affiliated with \(\M\) is
\(\varphi\)-measurable \cite[Proposition~4.7\textup{(iv)}]{hiai2021}.

For a general bounded input, the operator-valued series below are
interpreted as absolutely convergent series in \( L^1(\M,\varphi)\).
More precisely, for a series \(\sum_i x_i\) we prove that
\(\sum_i\|x_i\|_1<\infty\) and take its sum to be the resulting
\( L^1\)-limit.  This convergence also implies convergence in measure.
When the summands are positive, the finite partial sums increase, and this limit is also their supremum in the extended positive part of $\M$ by \cite[Proposition~4.23(2), Corollary~4.51 and Example~7.2(1)]{hiai2021} (see Section~\ref{subsec:extended-positive-part}).  Statements about these series use this convention unless another topology is specified.  Infinite sums of pairwise orthogonal projections, such as $\sum_ie_i=1$, are always understood in the strong operator topology.

\subsection{Center-valued dimension tools}
\label{subsec:center-valued-dimension-tools}

Let \(\mathcal N\) be a finite von Neumann algebra admitting a faithful
normal state, and let \(\Theta:\mathcal N\to Z(\mathcal N)\) be its normalized
center-valued trace.  The central type decomposition gives
pairwise orthogonal central projections \(z_{\mathrm{II}},z_1,z_2,\ldots\)
whose sum is \(1\) in the strong operator topology, where
\(z_m\mathcal N\) is homogeneous of type~\(\mathrm I_m\) and
\(z_{\mathrm{II}}\mathcal N\) is of type~\(\mathrm{II}_1\)
\cite[Sections~III.1.4.7 and~III.1.5.12--III.1.5.13]{MR2188261}.  The restriction of \(\Theta\)
to each nonzero central summand is its normalized center-valued trace.

For projections \(r,s\in\mathcal N\), we write \(r\sim s\) if there
is a partial isometry \(w\in\mathcal N\) with \(w^*w=r\) and
\(ww^*=s\).  We write \(r\precsim s\) if \(r\) is equivalent to a
subprojection of \(s\).  The comparison theorem for the center-valued
dimension function gives \(r\sim s\) if and only if \(\Theta(r)=\Theta(s)\)
and, more generally, \(r\precsim s\) if and only if
\(\Theta(r)\leq\Theta(s)\)
\cite[Theorem~8.4.3(vi), p.~532]{MR1468230}.
For \(\beta\in\Aut(\mathcal N)\) and \(x\in\mathcal N\), the uniqueness argument in the proof of
Proposition~\ref{prop:invariant-state-is-trace}, applied to
\(\beta^{-1}\circ\Theta\circ\beta\), gives
\(\Theta(\beta(x))=\beta(\Theta(x))\).

We will use these facts to choose subprojections with prescribed
center-valued traces.  On a homogeneous type-I summand, the possible
values are restricted by the matrix rank.

For a central positive contraction \(b\) in a homogeneous type~\(\mathrm I_m\) algebra \(\mathcal N\), the condition that \(mb\) is integer-valued means that there are pairwise orthogonal central projections
\(e_0,\ldots,e_m\) such that \(\sum_{j=0}^m e_j=1\) and
\(mb=\sum_{j=0}^m j e_j\). Equivalently, the spectrum of \(mb\) is
contained in \(\{0,1,\ldots,m\}\).

\begin{lemma}
\label{lem:type-I-dimension-range}
Suppose that \(\mathcal N\) is homogeneous of type~\(\mathrm I_m\).
Let \(r\in\Proj(\mathcal N)\) and \(b\in Z(\mathcal N)_+\) satisfy
\(0\leq b\leq\Theta(r)\).  There is a projection \(q\leq r\) with
\(\Theta(q)=b\) if and only if \(mb\) is integer-valued.
\end{lemma}

\begin{proof}
By the range theorem for the center-valued dimension function, a
central positive contraction \(b\) is the center-valued trace of a
projection in \(\mathcal N\) if and only if \(mb\) is integer-valued
\cite[Theorem~8.4.4(i), p.~533]{MR1468230}.  If a projection
\(q\leq r\) satisfies \(\Theta(q)=b\), this theorem applied to \(q\) gives
the required necessity.

Conversely, suppose that \(mb\) is integer-valued.  Since
\(0\leq b\leq\Theta(r)\leq\Theta(1)=1\), the range theorem provides a
projection \(s\in\mathcal N\) with \(\Theta(s)=b\).  The inequality
\(\Theta(s)\leq\Theta(r)\) and the comparison theorem give \(s\precsim r\)
\cite[Theorem~8.4.3(vi), p.~532]{MR1468230}.  Choose a partial isometry
\(w\in\mathcal N\) with \(w^*w=s\) and \(ww^*\leq r\), and let
\(q=ww^*\).  Then \(q\) is a subprojection of \(r\) equivalent to
\(s\), so traciality of \(\Theta\) gives
\(\Theta(q)=\Theta(s)=b\).
\end{proof}

\begin{lemma}
\label{lem:type-II-dimension-range}
Suppose that \(\mathcal N\) is a finite type-II von Neumann algebra
admitting a faithful normal state.  If \(r\in\Proj(\mathcal N)\) and
\(b\in Z(\mathcal N)_+\) satisfy \(0\leq b\leq\Theta(r)\), then there is
a projection \(q\leq r\) such that \(\Theta(q)=b\).
\end{lemma}

\begin{proof}
We have \(0\leq b\leq\Theta(r)\leq1\).  The type~\(\mathrm{II}_1\) range
theorem therefore gives a projection \(s\in\mathcal N\) with
\(\Theta(s)=b\) \cite[Theorem~8.4.4(ii), p.~533]{MR1468230}.  Since
\(\Theta(s)\leq\Theta(r)\), the comparison theorem gives \(s\precsim r\)
\cite[Theorem~8.4.3(vi), p.~532]{MR1468230}.  Thus there is a partial
isometry \(w\in\mathcal N\) with \(w^*w=s\) and \(ww^*\leq r\).
The projection \(q=ww^*\) lies below \(r\), and traciality of \(\Theta\)
gives \(\Theta(q)=\Theta(ww^*)=\Theta(w^*w)=\Theta(s)=b\).
\end{proof}

\begin{corollary}
\label{cor:center-valued-dimension-partition}
Let \(p\in\Proj(\mathcal N)\), and let \((c_i)_{i\geq1}\) be central
positive elements such that \(\sum_{i\geq1}c_i=\Theta(p)\) in the strong
operator topology.  Suppose that for every \(m\geq1\) and
\(i\geq1\), the element \(m z_m c_i\) is integer-valued in
\(Z(z_m\mathcal N)\). Then there are pairwise orthogonal projections
\(p_i\leq p\) such that for every \(i\geq1\), we have \(\Theta(p_i)=c_i\), and
\(\sum_{i\geq1}p_i=p\) strongly.
\end{corollary}

\begin{proof}
We construct the projections by induction.  Suppose that the mutually
orthogonal projections \(p_1,\ldots,p_{n-1}\leq p\) have been chosen
with \(\Theta(p_i)=c_i\), and let \(r_n=p-\sum_{i<n}p_i\).  For \(n=1\)
this means \(r_1=p\).  In either case \(r_n\) is a projection, and
linearity of \(\Theta\) and the hypothesis on the sum of the \(c_i\) give
\(\Theta(r_n)=\sum_{i\geq n}c_i\), with the sum taken strongly.  In
particular, \(0\leq c_n\leq\Theta(r_n)\).

For each nonzero type-I summand, let \(r_{n,m}=z_mr_n\) and
\(b_{n,m}=z_mc_n\).  Since \(\Theta\) is linear over the center,
\(\Theta(r_{n,m})=z_m\Theta(r_n)\), and hence
\(0\leq b_{n,m}\leq\Theta(r_{n,m})\).  The element \(mb_{n,m}\) is
integer-valued by assumption.  Applying
Lemma~\ref{lem:type-I-dimension-range} in \(z_m\mathcal N\), we obtain
a projection \(p_{n,m}\leq r_{n,m}\) with \(\Theta(p_{n,m})=b_{n,m}\).
Similarly, Lemma~\ref{lem:type-II-dimension-range}, applied in the
nonzero type-II summand, provides a projection
\(p_{n,\mathrm{II}}\leq z_{\mathrm{II}}r_n\) with
\(\Theta(p_{n,\mathrm{II}})=z_{\mathrm{II}}c_n\).  We take the corresponding
projection to be zero when a central summand is zero.

The projections \(p_{n,\mathrm{II}},p_{n,1},p_{n,2},\ldots\) lie in
mutually orthogonal central summands.  Their finite sums are therefore
projections bounded above by \(r_n\), so their strong sum
\(p_n=p_{n,\mathrm{II}}+\sum_{m\geq1}p_{n,m}\) is a projection below
\(r_n\).  By normality of \(\Theta\), its center-valued trace is
\(z_{\mathrm{II}}c_n+\sum_{m\geq1}z_mc_n=c_n\).  Moreover, \(p_n\)
is orthogonal to \(p_1,\ldots,p_{n-1}\), because it lies below
\(r_n=p-\sum_{i<n}p_i\).  This completes the induction.

The finite sums \(\sum_{i=1}^n p_i\) increase and are bounded above
by \(p\).  Let \(q=\sum_{i\geq1}p_i\) be their strong limit.  Then
\(q\) is a projection below \(p\), and normality gives
\(\Theta(q)=\sum_{i\geq1}\Theta(p_i)=\sum_{i\geq1}c_i=\Theta(p)\).  Thus
\(\Theta(p-q)=0\).  Since \(p-q\) is positive and \(\Theta\) is faithful, we
obtain \(p-q=0\), as required.
\end{proof}

\subsection{The extended positive part}
\label{subsec:extended-positive-part}

Let $\mathcal N$ be a von Neumann algebra, and let $\mathcal N_*^+$ be the cone of positive normal functionals on $\mathcal N$. Following Haagerup, in the formulation of \cite[Definition~7.1]{hiai2021}, the extended positive part $\widehat{\mathcal N}_+$ of $\mathcal N$ is the set of maps $m\colon\mathcal N_*^+\to[0,\infty]$ such that, for $\omega,\omega'\in\mathcal N_*^+$ and $\lambda\geq0$, we have $m(\omega+\omega')=m(\omega)+m(\omega')$ and $m(\lambda\omega)=\lambda m(\omega)$, and such that $m$ is lower semicontinuous for the norm topology on $\mathcal N_*^+$. Every $x\in\mathcal N_+$ defines the element $\omega\mapsto\omega(x)$ of $\widehat{\mathcal N}_+$, and in this way we regard $\mathcal N_+$ as a subset of $\widehat{\mathcal N}_+$. For $m,m'\in\widehat{\mathcal N}_+$, we let $(m+m')(\omega):=m(\omega)+m'(\omega)$, and we write $m\leq m'$ if $m(\omega)\leq m'(\omega)$ for every $\omega\in\mathcal N_*^+$. If $(m_i)_i$ is an increasing net in $\widehat{\mathcal N}_+$, then its supremum belongs to $\widehat{\mathcal N}_+$ and is given by $(\sup_i m_i)(\omega)=\sup_i m_i(\omega)$ \cite[Section~7.1, p.~95]{hiai2021}.

Let $\Psi\colon\mathcal N\to\mathcal N'$ be a normal positive linear map between von Neumann algebras. For $m\in\widehat{\mathcal N}_+$ and $\omega\in\mathcal N'^{~+}_{*}$, define $\widehat\Psi(m)(\omega):=m(\omega\circ\Psi)$. Since the map $\omega\mapsto\omega\circ\Psi$ is positive, linear and norm continuous, $\widehat\Psi(m)$ belongs to $\widehat{\mathcal N'}_+$ by \cite[Definition~7.1]{hiai2021}.

Suppose now that $\mathcal N$ is finite with a faithful normal finite trace $\tau$. For $m\in\widehat{\mathcal N}_+$, we let $\widehat\tau(m):=m(\tau)\in[0,\infty]$. For $h\in L^1(\mathcal N,\tau)_+$, we write $\tau(h~\cdot~)$ for the positive normal functional $x\mapsto\tau(hx)$ \cite[Corollary~4.51]{hiai2021}. For $y\in L^1(\mathcal N,\tau)_+$ and $k\in\N$, let $y\wedge k:=f_k(y)$, where $f_k(t)=\min\{t,k\}$, and let $\widehat y:=\sup_{k\geq1}(y\wedge k)\in\widehat{\mathcal N}_+$, as in \cite[Example~7.2(1)]{hiai2021}.

\section{Non-commutative Kac's lemma for \texorpdfstring{$\mathbb Z$-action}{}}
\label{sec:Z-action}

In this section we study the non-commutative analogue of Kac's lemma for actions of the countable abelian group $\Z$. Later we extend the results to actions of countable discrete groups. Let us fix a finite von Neumann algebra $\M$ equipped with a faithful, normal (f.n.) state $\varphi$. Let $\Gamma$ denote a group and $\alpha: \Gamma \to \Aut(\M)$ be an action. We further assume that the action $\alpha$ preserves the underlying state $\varphi$. We assume that $\Gamma$ acts ergodically on $\M$. By Proposition \ref{prop:invariant-state-is-trace}, $\varphi$ is tracial. Thus throughout this section, $(\M,\varphi)$ is a finite von Neumann algebra with a f.n. tracial state, and $\alpha$ is an ergodic $\varphi$-preserving action.

\subsection{Action of the group \texorpdfstring{$\Z$}{}:} Let us assume that $\Gamma=\Z$ and fix a non-zero projection $p \in \M$. Let \(q_0=s_0=0\).  For \(m\in \mathbb{N}\), define
\begin{align}
q_m
  &:=\bigvee_{k=1}^{m}\alpha^{-k}(p),
&
r_m
  &:=q_m-q_{m-1},                                      \label{eq:Z-QR}\\
s_m
  &:=\bigvee_{k=1}^{m}\alpha^k(p),
&
\widetilde t_m
  &:=(s_{m-1}\vee p)-s_{m-1},                          \label{eq:Z-Stilde}\\
t_m
  &:=p-(p\wedge s_{m-1}),
&
f_m
  &:=t_m-t_{m+1}.                                      \label{eq:Z-TF}
\end{align}
Observe that since $\alpha$ is an action by $*$-automorphisms on the von Neumann algebra $\M$, $\alpha^k(p)$ is a projection for all $k \in \Z$ and so $q_m$ and $s_m$ are projections for all $m \in \N$. Since $(q_m)_m$ is an increasing sequence of projections, $r_m$ is a projection for every $m\in\N$. Similarly, since $(t_m)_m$ is a decreasing sequence of projections, $f_m$ is a projection in $\M$ for every $m\in\N$.

\begin{proposition}
\label{prop:Z-projection-structure}
The projections \(q_m,r_m,s_m,\widetilde t_m,t_m,f_m\)
defined in \eqref{eq:Z-QR}--\eqref{eq:Z-TF} have the following
properties.
\begin{enumerate}
\item For $m\neq n$, we have $r_m r_n= r_n r_m=0$, and $\sum_{m=1}^{\infty}r_m=1$.
\item For every $m\geq1$, we have $\alpha^m(r_m)=\widetilde t_m$.
\item There is a partial isometry $u_m\in\M$ such that $u_m^*u_m=t_m$ and $u_mu_m^*=\widetilde t_m$.
\item For $m\neq n$, we have $f_m f_n= f_n f_m=0$. Moreover, $t_m=\sum_{n=m}^{\infty}f_n$ and $p=\sum_{n=1}^{\infty}f_n$.
\item For every $m$, we have $\varphi(t_m)=\varphi(\widetilde t_m)=\varphi(r_m)$, and consequently $\sum_{m=1}^{\infty}\varphi(t_m)=1$.
\end{enumerate}
\end{proposition}

\begin{proof}
The sequence \((q_m)_m\) is increasing, so \(r_m\) are
pairwise orthogonal.  Let $q_\infty:=\bigvee_{k\geq1}\alpha^{-k}(p)$.
Then \(\alpha^{-1}(q_\infty)\leq q_\infty\), while state invariance
gives
\(\varphi(\alpha^{-1}(q_\infty))=\varphi(q_\infty)\). Therefore, faithfulness of $\varphi$ implies
\(\alpha^{-1}(q_\infty)=q_\infty\).  The projection \(q_\infty\) is
nonzero and \(\alpha\)-invariant, so ergodicity gives \(q_\infty=1\).
Since
\(\sum_{m=1}^{N}r_m=q_N\), the first assertion follows.

Now, for \(m\geq1\), observe that
\begin{align*}
\alpha^m(r_m)
 &=\alpha^m(q_m)-\alpha^m(q_{m-1})\\
 &=\bigvee_{k=0}^{m-1}\alpha^k(p)
   -\bigvee_{k=1}^{m-1}\alpha^k(p)\\
 &=(s_{m-1}\vee p)-s_{m-1}
 =\widetilde t_m.
\end{align*}
By the parallelogram law
\cite[Proposition~III.1.1.3]{MR2188261}, we obtain the required partial isometry $u_m \in \M$. Therefore, traciality and state invariance give $\varphi(t_m)=\varphi(\widetilde t_m)=\varphi(\alpha^m(r_m))=\varphi(r_m)$.
In particular, \(\varphi(t_m)\to0\) as $m \to \infty$.  If
\(t_\infty=\bigwedge_mt_m\), normality and faithfulness of $\varphi$ imply \(t_\infty=0\).
Therefore, $\sum_{n=m}^{N}f_n=t_m-t_{N+1}\longrightarrow t_m$ as $N\to\infty$. Now taking \(m=1\), and observing that \(t_1=p\), proves the
fourth assertion.  Finally, $\sum_{m\geq1}\varphi(t_m)=\sum_{m\geq1}\varphi(r_m)=1$.
\end{proof}

Fix the partial isometries $(u_m)_m$ and the projections from Proposition~\ref{prop:Z-projection-structure}. Since $\M$ is finite and $\varphi$ is tracial, we have $\M\subseteq L^1(\M,\varphi)$. We use these projections and partial isometries to form the following series along the orbit of $x\in\M$.

For $x \in \M$, let
$\widetilde{K}_p(x):=\sum_{m=1}^{\infty}\widetilde t_m\alpha^m(x)\widetilde t_m$ and $K_{p,u}(x):=\sum_{m=1}^{\infty}u_m^*\alpha^m(x)u_m$.

\begin{proposition}
\label{thm:Z-Kac-transport}
For \(x\in\M\), the following hold true. We have
$\sum_{m=1}^{\infty}\|\widetilde t_m\alpha^m(x)\widetilde t_m\|_1\leq\|x\|$ and $\sum_{m=1}^{\infty}\|u_m^*\alpha^m(x)u_m\|_1\leq\|x\|$,
where $\|x\|$ denotes the norm in $\M$. Furthermore, we have $\varphi(\widetilde{K}_p(x))=\varphi(K_{p,u}(x))=\varphi(x)$.
\end{proposition}

\begin{proof}
We first note that for all $m \in \N$, $|\widetilde t_m \alpha^m(x) \widetilde t_m| \leq \|x\| \widetilde t_m$. Thus $\|\widetilde t_m \alpha^m(x) \widetilde t_m\|_1 \leq \|x\| \varphi(\widetilde t_m)$. Similarly, since $ t_m$ is the initial projection of the partial isometry $u_m$, we have $|u_m^*\alpha^m(x)u_m|\leq\|x\|t_m$ and hence $\|u_m^*\alpha^m(x)u_m\|_1\leq\|x\|\varphi(t_m)$. Since $\sum_m\varphi(\widetilde t_m)=\sum_m\varphi(t_m)=1$, we have the required norm inequalities. Moreover, for each \(m\),
\(\alpha^m(r_m)=\widetilde t_m\) implies
\begin{align*}
\varphi\bigl(
  \widetilde t_m\alpha^m(x)\widetilde t_m
\bigr)=\varphi\bigl(\alpha^m(r_mxr_m)\bigr)=\varphi(r_mxr_m)
 =\varphi(xr_m).
\end{align*}
Similarly,
$\varphi\bigl(u_m^*\alpha^m(x)u_m\bigr)
=\varphi\bigl(\alpha^m(x)u_mu_m^*\bigr)=\varphi\bigl(\alpha^m(x)\widetilde t_m\bigr)
 =\varphi(xr_m)$.
Since \(q_N=\sum_{m=1}^{N}r_m\) converges to \(1\), it follows that
$\sum_{m=1}^{N}\varphi(xr_m)
  =\varphi(xq_N)\longrightarrow\varphi(x)$.
Continuity of \(\varphi\) on \(L^1\) gives the asserted trace identities.
\end{proof}

We now use these projections to prove the non-commutative analogue of Kac's formula \cite{kac1947}. Let $T_p:=\sum_{n=1}^{\infty}nf_n$, and let $\varphi_p(\cdot):=\varphi(\cdot)/\varphi(p)$ be the normalized state on $p\M p$.

\begin{theorem}[Noncommutative Kac formula]\label{thm:quantum-return-time}
	The operator $T_p$
	is a positive self-adjoint $\varphi$-measurable operator affiliated with
	$p\mathcal{M}p$, belongs to $ L^1(p\M p)_+$, and satisfies $T_p=\sum_{m=1}^{\infty}t_m$. Furthermore, we have
	\(
	\varphi(T_p)=1.
	\)
	Therefore, $\varphi_p(T_p)=1/\varphi(p)$.
\end{theorem}

\begin{proof}
	By Proposition~\ref{prop:Z-projection-structure}, the projections  \(f_n\) are	pairwise orthogonal and $\sum_{n \in \N} f_n =p$.  The spectral theorem
	therefore defines the positive affiliated operator
	\(T_p=\sum_n nf_n\), with $f_n=\one_{\{n\}}( T_p)$ and $t_m=\one_{[m,\infty)}( T_p)$. Therefore, for any $N\geq1$, we have $\sum_{m=1}^{N}t_m=\sum_{n=1}^{\infty}\min\{n,N\}f_n$.
	The right-hand side increases to \(\sum_n nf_n\).  On applying
	\(\varphi\) and using
	Proposition~\ref{prop:Z-projection-structure}, we obtain
	$\varphi(T_p)=\sum_{m=1}^{\infty}\varphi(t_m)=\sum_{m=1}^{\infty}\varphi(r_m)=1$.
	Thus by Proposition \ref{thm:Z-Kac-transport} the positive series converges in \( L^1\), and
	\(K_{p,u}(1)=\sum_m u_m^*u_m=\sum_mt_m= T_p\).
	The normalized formula follows by division by \(\varphi(p)\).
\end{proof}
For $m\geq0$, let $\mathcal A_m:=W^*(p,\alpha(p),\ldots,\alpha^m(p))$ and $\mathcal N_m:=p\mathcal A_m p$.
Let $e_0:=0$, and for $m\geq1$, define $e_m:=p\wedge s_m=p-t_{m+1}$.

\begin{corollary}
\label{cor:adapted-return}
For $m\geq0$, we have $e_m\in\Proj(\mathcal N_m)$ and $e_m\leq e_{m+1}$. Moreover, $e_m\longrightarrow p$, and for $m\geq1$, we have $f_m=e_m-e_{m-1}$.
\end{corollary}

\begin{proof}
The projection \(p\wedge s_m\) belongs to the von Neumann algebra
\(\mathcal A_m\) and is dominated by \(p\), hence belongs to the corner
\(\mathcal N_m\).  The remaining assertions follow from the
monotonicity of \(s_m\), the identity \(e_m=p-t_{m+1}\), and
\(t_m\downarrow0\).
\end{proof}

\begin{remark}
\label{rem:stopping-time-scope}
Observe that the family $(e_n)_n$ can be seen as the stopping time associated to the von Neumann subalgebra filtration $(\mathcal N_n)_n$. Accordingly, the word
``quantum'' refers here to the replacement of measurable events by
projections and of an integer-valued random variable by the positive
affiliated operator \(T_p\) with its spectral resolution.  No
repeated-measurement dynamics is implicit, and
\(\varphi_p(f_n)\) is not asserted to be the first-detection
probability of a measurement protocol.  Such an interpretation would
require a separate update rule and a proof that it produces the same
projections.  The identity above is likewise distinct from averaged
recurrence statements such as those in \cite{nsz03}, from
noncommutative stopping-time and time-projection frameworks
\cite{barnett-thakrar,luczak2005}, and from recurrence times for
iterated open quantum systems \cite{sinkovicz2016}.
\end{remark}

In the following proposition, we show that the non-commutative Kac's lemma described above reduces to the classical Kac formula. Let $T$ be an invertible ergodic measure-preserving transformation of a probability space $(X,\mathcal B,\mu)$, and let $A\in\mathcal B$ satisfy $\mu(A)>0$. Let $\M=L^\infty(X,\mu)$, $\varphi(f)=\int_X f~d\mu$, $\alpha(f)=f\circ T$ and $p=\one_A$. For $z\in A$, define $r_A(z):=\inf\{m\geq1:T^m z\in A\}$, with $\inf\varnothing=\infty$.
\begin{proposition}
	\label{prop:commutative-Z}
	For the projections $t_m,f_m$ defined in \eqref{eq:Z-TF} and the operator $T_p$ of Theorem~\ref{thm:quantum-return-time}, we have, modulo null sets, $t_m=\one_{\{z\in A:r_A(z)\geq m\}}$, $f_m=\one_{\{z\in A:r_A(z)=m\}}$, and $T_p=r_A$.
	For the canonical choice $u_m=t_m=\widetilde t_m$ and $f\in L^\infty(X,\mu)$, we have, for almost every $z\in A$, $K_{p,u}(f)(z)=\sum_{m=1}^{r_A(z)}f(T^m z)$,
	and $\int_A\sum_{m=1}^{r_A(z)}f(T^m z)~d\mu(z)=\int_X f~d\mu$.
	In particular, $\int_A r_A~d\mu=1$, which is the classical Kac formula \cite{kac1947}.
\end{proposition}

\begin{proof}
	The support of \(s_{m-1}\) is
	\(\bigcup_{k=1}^{m-1}T^{-k}(A)\).  Therefore, $t_m=\one_{A\setminus\bigcup_{k=1}^{m-1}T^{-k}(A)}=\one_{\{r_A\geq m\}}$.
	Taking consecutive differences gives the formula for \(f_m\), and
	Theorem~\ref{thm:quantum-return-time} gives \(T_p=r_A\). Furthermore, \(\widetilde t_m=t_m\), and \(u_m=t_m\) is an admissible choice. Since $K_{p,u}(x)=
  \sum_{m=1}^{\infty}
  u_m^*\alpha^m(x)u_m$, the orbit-sum identity follows; the integral identity follows from
	Proposition~\ref{thm:Z-Kac-transport}.
\end{proof}

\section{Non-commutative Kac's lemma for general countable group action}
\label{sec:universal-towers}

In this section, we extend Kac's lemma to actions of countable groups. Let $\Gamma$ be a countable discrete group with identity $e$, and let $\alpha\colon\Gamma\to\Aut(\M)$ be an ergodic $\varphi$-preserving action. By Proposition~\ref{prop:invariant-state-is-trace}, $\varphi$ is tracial, so its centralizer $\M_\varphi$ is all of $\M$. Fix a non-zero projection $p\in\M$. Observe that $\bigvee_{\gamma\in\Gamma}\alpha_\gamma(p)=1$.
Indeed, let $q_\infty:=\bigvee_{\gamma\in\Gamma}\alpha_\gamma(p)$. For $\eta\in\Gamma$, the automorphism $\alpha_\eta$ preserves suprema of projections, and $\gamma\mapsto\eta\gamma$ is a bijection of $\Gamma$; hence $\alpha_\eta(q_\infty)=\bigvee_{\gamma\in\Gamma}\alpha_{\eta\gamma}(p)=q_\infty$.
Thus $q_\infty\in\M^\alpha=\C1$, and since $q_\infty\geq p\neq0$, we obtain $q_\infty=1$.

Throughout this section, we let $J:=\N$ if $\Gamma$ is infinite and $J:=\{1,\ldots,N\}$ if $|\Gamma|=N<\infty$, and we fix a bijective enumeration $\Gamma=\{g_j:j\in J\}$ with $g_1=e$.

\subsection{Weak partial-isometric allocations}

Let $T\colon\Gamma\curvearrowright(X,\mathcal B,\mu)$ be a measure-preserving action on a probability space, and let $A\in\mathcal B$. Meyerovitch and Weiss call a measurable map $\kappa\colon X\to\Gamma$ an $A$-allocation if $T_{\kappa(y)}(y)\in A$ for almost every $y\in X$ \cite[Definition~1.2]{mw24}. In terms of the action $\alpha_\gamma(f)=f\circ T_\gamma^{-1}$ on $L^\infty(X,\mu)$, the projections $\rho_\gamma:=\one_{\kappa^{-1}(\gamma)}$ are pairwise orthogonal, $\sum_\gamma\rho_\gamma=1$, and $\alpha_\gamma(\rho_\gamma)\leq\one_A$ for every $\gamma\in\Gamma$ (Proposition~\ref{prop:MW-comparison}). For a non-commutative algebra $\M$, projections with the analogous properties need not exist (Example~\ref{ex:direct-allocation-fails}), and we therefore only require $\alpha_\gamma(\rho_\gamma)$ to be equivalent to a subprojection of $p$.

\begin{definition}
\label{def:partial-isometric-allocation}
A family $(\rho_\gamma)_{\gamma\in\Gamma}$ of pairwise orthogonal projections in $\M$ with $\sum_{\gamma\in\Gamma}\rho_\gamma=1$ is called a weak partial-isometric $p$-allocation, or simply a partial-isometric $p$-allocation, for $\alpha$ if, for every $\gamma\in\Gamma$, there are projections $p_\gamma\leq p$ and partial isometries $u_\gamma\in\M$ such that $u_\gamma^*u_\gamma=p_\gamma$ and $u_\gamma u_\gamma^*=\alpha_\gamma(\rho_\gamma)$.
We call $(p_\gamma,u_\gamma)_{\gamma\in\Gamma}$ a transport witness for the allocation. If, moreover, $\alpha_\gamma(\rho_\gamma)\leq p$ for every $\gamma\in\Gamma$, we call $(\rho_\gamma)_{\gamma\in\Gamma}$ a direct $p$-allocation; in this case, $p_\gamma=u_\gamma:=\alpha_\gamma(\rho_\gamma)$ defines a transport witness, so every direct $p$-allocation is a partial-isometric $p$-allocation. For a partial-isometric $p$-allocation with transport witness $(p_\gamma,u_\gamma)_{\gamma\in\Gamma}$ and $x\in\M$, we consider the series $L_{\rho,p,u}(x):=\sum_{\gamma\in\Gamma}u_\gamma^*\alpha_\gamma(x)u_\gamma$.
\end{definition}

For the following theorem, $(\M,\varphi)$ is a finite von Neumann algebra with a faithful normal tracial state, and $\alpha\colon\Gamma\to\Aut(\M)$ is an action of a countable discrete group preserving $\varphi$. Ergodicity is not assumed.

\begin{theorem}
\label{thm:order-free-allocation}
Let $(\rho_\gamma)_{\gamma\in\Gamma}$ be a partial-isometric $p$-allocation with transport witness $(p_\gamma,u_\gamma)_{\gamma\in\Gamma}$. For every $x\in\M$, we have $\sum_{\gamma\in\Gamma}\|u_\gamma^*\alpha_\gamma(x)u_\gamma\|_1\leq\|x\|$.
In particular, the series $L_{\rho,p,u}(x)$ converges absolutely in $L^1(p\M p,\varphi)$. Furthermore, $\varphi(L_{\rho,p,u}(x))=\varphi(x)$.
Moreover, if $x\geq0$, then $L_{\rho,p,u}(x)\in L^1(p\M p,\varphi)_+$, and $\widehat{L_{\rho,p,u}(x)}$ is the supremum of the increasing net of finite partial sums in the extended positive part of $p\M p$ (see Section~\ref{subsec:extended-positive-part}).
\end{theorem}

\begin{proof}
For $\gamma\in\Gamma$, let $a_\gamma:=u_\gamma^*\alpha_\gamma(x)u_\gamma$. Since $u_\gamma=u_\gamma u_\gamma^*u_\gamma=u_\gamma p_\gamma$, we have $a_\gamma=p_\gamma a_\gamma p_\gamma\in p_\gamma\M p_\gamma$. Also, $\|a_\gamma\|\leq\|u_\gamma\|^2\|\alpha_\gamma(x)\|\leq\|x\|$, because the $*$-automorphism $\alpha_\gamma$ is isometric. The positive element $|a_\gamma|$ of the von Neumann algebra $p_\gamma\M p_\gamma$, whose unit is $p_\gamma$, satisfies $|a_\gamma|\leq\|a_\gamma\|p_\gamma$, and hence
\begin{equation}
\label{eq:corner-trace-estimate}
\|a_\gamma\|_1=\varphi(|a_\gamma|)\leq\|a_\gamma\|~\varphi(p_\gamma)\leq\|x\|~\varphi(p_\gamma).
\end{equation}
By traciality and invariance of $\varphi$,
\[
\varphi(p_\gamma)=\varphi(u_\gamma^*u_\gamma)=\varphi(u_\gamma u_\gamma^*)=\varphi\bigl(\alpha_\gamma(\rho_\gamma)\bigr)=\varphi(\rho_\gamma),
\]
and since $\varphi$ is normal and $\sum_\gamma\rho_\gamma=1$, we obtain $\sum_\gamma\varphi(p_\gamma)=\sum_\gamma\varphi(\rho_\gamma)=1$. Hence $\sum_\gamma\|a_\gamma\|_1\leq\|x\|$.

For a finite set $F\subseteq\Gamma$, let $S_F:=\sum_{\gamma\in F}a_\gamma$ and $\rho_F:=\sum_{\gamma\in F}\rho_\gamma$. Let $\varepsilon>0$, and choose a finite set $F_0\subseteq\Gamma$ with $\sum_{\gamma\notin F_0}\|a_\gamma\|_1<\varepsilon$. For finite sets $F,F'\supseteq F_0$, we have $\|S_F-S_{F'}\|_1\leq\sum_{\gamma\in F\triangle F'}\|a_\gamma\|_1<\varepsilon$.
Thus $(S_F)_F$ is a Cauchy net in the Banach space $L^1(p\M p,\varphi)$, and it converges to an element $L_{\rho,p,u}(x)$. In particular, the partial sums along any enumeration of $\Gamma$ converge to the same limit.

Using traciality, the identity $u_\gamma u_\gamma^*=\alpha_\gamma(\rho_\gamma)$ and invariance of $\varphi$, we obtain
\[
\varphi(S_F)=\sum_{\gamma\in F}\varphi\bigl(\alpha_\gamma(x)u_\gamma u_\gamma^*\bigr)=\sum_{\gamma\in F}\varphi\bigl(\alpha_\gamma(x\rho_\gamma)\bigr)=\varphi(x\rho_F).
\]
Let $e_F:=1-\rho_F$. Traciality gives $\varphi(x)-\varphi(x\rho_F)=\varphi(xe_F)=\varphi(e_Fxe_F)$, and the argument leading to \eqref{eq:corner-trace-estimate}, applied to $e_Fxe_F\in e_F\M e_F$ in place of $a_\gamma\in p_\gamma\M p_\gamma$, gives $|\varphi(e_Fxe_F)|\leq\|e_Fxe_F\|_1\leq\|x\|~\varphi(e_F)$.
By normality, $\varphi(e_F)\to0$. For $y\in L^1(p\M p,\varphi)$, we have $|\varphi(y)|\leq\|y\|_1$, so we conclude that $\varphi(L_{\rho,p,u}(x))=\lim_F\varphi(S_F)=\lim_F\varphi(x\rho_F)=\varphi(x)$.

Finally, let $x\geq0$. Then every $a_\gamma$ is positive, so $(S_F)_F$ is an increasing net in $(p\M p)_+$ which converges in $L^1(p\M p,\varphi)$ to $L_{\rho,p,u}(x)$. The identification with the positive predual in \cite[Corollary~4.51]{hiai2021} shows that the positive cone is norm closed, so $L_{\rho,p,u}(x)\geq0$. Since $\Gamma$ is countable, the net $(S_F)_F$ has a cofinal sequence. For $h\in L^1(p\M p,\varphi)_+$, the spectral identification \cite[Example~7.2(1)]{hiai2021} and monotone convergence \cite[Proposition~4.23(2) and Corollary~4.21(2)]{hiai2021} give
\begin{align*}
\sup_F S_F(\varphi(h~\cdot~))
&=\sup_F\varphi(h^{1/2}S_Fh^{1/2})\\
&=\varphi\bigl(h^{1/2}L_{\rho,p,u}(x)h^{1/2}\bigr)
=\widehat{L_{\rho,p,u}(x)}(\varphi(h~\cdot~)).
\end{align*}
Every positive normal functional on $p\M p$ has this form by \cite[Corollary~4.51]{hiai2021}, so $\widehat{L_{\rho,p,u}(x)}=\sup_FS_F$.
\end{proof}

We use Theorem~\ref{thm:order-free-allocation} for actions that need not be ergodic in Proposition~\ref{prop:MW-comparison} and Remark~\ref{rem:towers-allocation}.

\begin{remark}
\label{rem:Z-case-allocation}
The assertions about $K_{p,u}$ in Proposition~\ref{thm:Z-Kac-transport} are the case $\Gamma=\Z$ of Theorem~\ref{thm:order-free-allocation}. Indeed, let $\alpha$ be as in Section~\ref{sec:Z-action}, so that for $m\in\Z$, we have $\alpha_m=\alpha^m$, and let $r_m$, $s_m$, $t_m$, $\widetilde t_m$ and $u_m$ be as in \eqref{eq:Z-QR}--\eqref{eq:Z-TF} and Proposition~\ref{prop:Z-projection-structure}. For $m\geq1$, let $\rho_m:=r_m$ and $p_m:=t_m$, and for $m\leq0$, let $\rho_m:=p_m:=u_m:=0$. By Proposition~\ref{prop:Z-projection-structure}, the projections $\rho_m$ are pairwise orthogonal with $\sum_{m\in\Z}\rho_m=1$, and for $m\geq1$, we have $p_m=p-(p\wedge s_{m-1})\leq p$, $u_m^*u_m=t_m=p_m$ and $u_mu_m^*=\widetilde t_m=\alpha^m(r_m)=\alpha_m(\rho_m)$.
Hence $(\rho_m)_{m\in\Z}$ is a partial-isometric $p$-allocation with transport witness $(p_m,u_m)_{m\in\Z}$. For $x\in\M$, we have $L_{\rho,p,u}(x)=K_{p,u}(x)$. In particular, $L_{\rho,p,u}(1)=\sum_{m\geq1}t_m=T_p$, and Theorem~\ref{thm:order-free-allocation} with $x=1$ gives the identity $\varphi(T_p)=1$ of Theorem~\ref{thm:quantum-return-time}.
\end{remark}

\begin{remark}
\label{rem:allocation-identity-scope}
By Theorem~\ref{thm:order-free-allocation}, the value $\varphi(L_{\rho,p,u}(x))=\varphi(x)$ does not depend on the allocation and the transport witness, whereas the operator $L_{\rho,p,u}(x)$ does. For instance, let $\Z$ act on $X=\Z/3\Z$, equipped with the uniform probability measure, by $T_\gamma(y)=y+\gamma$; this action is ergodic, since it is transitive. Let $\M=L^\infty(X)$, $\alpha_\gamma(f)=f\circ T_\gamma^{-1}$ and $p=\one_A$ with $A=\{0,1\}$, and note that $\alpha_\gamma(\one_E)=\one_{T_\gamma(E)}$ for $E\subseteq X$. The families given by $\rho_0=\one_{\{0,1\}}$, $\rho_1=\one_{\{2\}}$, and by $\rho_0=\one_{\{0,1\}}$, $\rho_2=\one_{\{2\}}$, with all other $\rho_\gamma$ equal to zero, are direct $p$-allocations, since $\alpha_1(\one_{\{2\}})=\one_{\{0\}}$ and $\alpha_2(\one_{\{2\}})=\one_{\{1\}}$. With the transport witnesses $p_\gamma=u_\gamma=\alpha_\gamma(\rho_\gamma)$, we obtain
\[
L_{\rho,p,u}(1)=\sum_{\gamma\in\Z}p_\gamma=2\cdot\one_{\{0\}}+\one_{\{1\}}~\text{and}~L_{\rho,p,u}(1)=\one_{\{0\}}+2\cdot\one_{\{1\}},
\]
respectively. For $\Gamma=\Z$, Remark~\ref{rem:Z-case-allocation} shows that the allocation constructed in Section~\ref{sec:Z-action} from the positive powers of $\alpha$ gives $L_{\rho,p,u}(1)=T_p$. A general countable group has no distinguished positive-time order that would single out one allocation, and Theorem~\ref{thm:order-free-allocation} holds for each of them.
\end{remark}

We now use the translates of $p$ to construct orthogonal source projections for a partial-isometric $p$-allocation.

\begin{lemma}
\label{lem:greedy-quantum-allocation}
Let $\alpha\colon\Gamma\to\Aut(\M)$ be an action of $\Gamma$ on a von Neumann algebra $\M$, not necessarily ergodic, and let $p\in\Proj(\M)$ satisfy $\bigvee_{\gamma\in\Gamma}\alpha_\gamma(p)=1$. Then there are pairwise orthogonal projections $(e_j)_{j\in J}$ with $\sum_{j\in J}e_j=1$, and projections $d_j\leq\alpha_{g_j}^{-1}(p)$, such that for every $j\in J$, we have $e_j\sim d_j$, and $e_1=d_1=p$.
\end{lemma}

\begin{proof}
For $j\in J$, let $a_j:=\alpha_{g_j}^{-1}(p)$ and $k_j:=\bigvee_{i=1}^ja_i$, and let $k_0:=0$. For $j\in J$, consider the projections $e_j:=k_j-k_{j-1}$.
Since $(k_j)_j$ is increasing, the $e_j$ are pairwise orthogonal and $\sum_{i\leq j}e_i=k_j$. As $\alpha_{g_j}^{-1}=\alpha_{g_j^{-1}}$ and $j\mapsto g_j^{-1}$ is a bijection from $J$ onto $\Gamma$, we have $\bigvee_{j\in J}k_j=\bigvee_{\gamma\in\Gamma}\alpha_\gamma(p)=1$. An increasing sequence of projections converges strongly to its supremum, and hence $\sum_{j\in J}e_j=1$. Since $k_j=a_j\vee k_{j-1}$, the parallelogram law \cite[Proposition~III.1.1.3]{MR2188261} gives $e_j=(a_j\vee k_{j-1})-k_{j-1}\sim a_j-(a_j\wedge k_{j-1})=:d_j$, and clearly $d_j\leq a_j$. Finally, $g_1=e$ gives $a_1=k_1=p$, so $e_1=p$ and $d_1=p-(p\wedge0)=p$. Applying $\alpha_{g_j}$ to $d_j\leq a_j$ gives $\alpha_{g_j}(d_j)\leq p$.
\end{proof}

For the projections $d_j$ supplied by Lemma~\ref{lem:greedy-quantum-allocation}, let $p_j:=\alpha_{g_j}(d_j)$. Then $p_j\leq p$ for every $j\in J$, and $p_1=p$.

\begin{corollary}
\label{cor:allocation-existence}
For an ergodic $\varphi$-preserving action of a countable discrete group on a finite von Neumann algebra $\M$ with a faithful normal state $\varphi$, every non-zero projection $p\in\M$ admits a partial-isometric $p$-allocation.
\end{corollary}

\begin{proof}
Since $\bigvee_{\gamma\in\Gamma}\alpha_\gamma(p)=1$, Lemma~\ref{lem:greedy-quantum-allocation} applies to $p$.  For each $j$, choose a
partial isometry $s_j\in\M$ with $s_j^*s_j=e_j$ and
$s_js_j^*=d_j$.  Let $\rho_{g_j}=e_j$, $p_{g_j}=p_j$ and
$u_{g_j}=\alpha_{g_j}(s_j^*)$.  Then
\begin{align*}
  u_{g_j}^*u_{g_j}
  &=\alpha_{g_j}(s_js_j^*)
    =\alpha_{g_j}(d_j)=p_j\leq p,\\
  u_{g_j}u_{g_j}^*
  &=\alpha_{g_j}(s_j^*s_j)=\alpha_{g_j}(e_j)
    =\alpha_{g_j}(\rho_{g_j}).
\end{align*}
The projections $\rho_{g_j}$ are pairwise orthogonal and have
sum $1$, so these families satisfy
Definition~\ref{def:partial-isometric-allocation}.
\end{proof}

\begin{remark}
\label{rem:direct-allocation}
For non-commutative algebras, a partial-isometric $p$-allocation need not be direct, and a direct $p$-allocation need not exist, as the following example shows. If $\M$ is abelian, the two notions of allocation coincide. Indeed, a partial isometry $u$ in an abelian von Neumann algebra satisfies $u^*u=uu^*$, so equivalent projections in $\M$ are equal, and a transport witness of a partial-isometric $p$-allocation satisfies $\alpha_\gamma(\rho_\gamma)=u_\gamma u_\gamma^*=u_\gamma^*u_\gamma=p_\gamma\leq p$.
\end{remark}

\begin{example}
\label{ex:direct-allocation-fails}
Let $\Gamma=\Z_2^2$ act on $M_2(\C)$ by conjugation with the
four Pauli matrices
\[
  I,~
  \sigma_x=\begin{pmatrix}0&1\\1&0\end{pmatrix},~
  \sigma_y=\begin{pmatrix}0&-\mathrm i\\ \mathrm i&0\end{pmatrix},~
  \sigma_z=\begin{pmatrix}1&0\\0&-1\end{pmatrix}.
\]
Note that each $\sigma_i$ satisfies $\sigma_i^2=I$, and the product of two distinct
Pauli matrices is a scalar multiple of the third.  Since scalar
multiples give the same conjugation, these four conjugations define an
action of $\Z_2^2$.  A matrix fixed by conjugation with
$\sigma_z$ is diagonal.  If it is also fixed by conjugation with
$\sigma_x$, its two diagonal entries are equal.  The fixed-point
algebra is therefore $\C I$, and the action preserves the
normalized trace.

Let $v=(1,1,1)/\sqrt3$, and let
$p=\frac12(I+v_1\sigma_x+v_2\sigma_y+v_3\sigma_z)$.
The Pauli relations give
$(v_1\sigma_x+v_2\sigma_y+v_3\sigma_z)^2=I$, so $p$ is a
projection.  Its unnormalized matrix trace is $1$, and hence it
has rank one.  Conjugating by $\sigma_i$ preserves its own
coefficient and reverses the other two.  Thus the four orbit
projections have Bloch vectors
\[
  \frac1{\sqrt3}(1,1,1),~
  \frac1{\sqrt3}(1,-1,-1),~
  \frac1{\sqrt3}(-1,1,-1),~
  \frac1{\sqrt3}(-1,-1,1).
\]
Any two distinct vectors $w,w'$ in this list satisfy
$w\cdot w'=-1/3$.  The identities
$\Tr(\sigma_i)=0$ and $\Tr(\sigma_i\sigma_j)=2\delta_{ij}$
therefore give $\Tr(p_wp_{w'})=(1+w\cdot w')/2=1/3$ for the corresponding
orbit projections $p_w=\frac12(I+w_1\sigma_x+w_2\sigma_y+w_3\sigma_z)$ and $p_{w'}$.  In particular, no two distinct orbit
projections are orthogonal.

Suppose that a direct $p$-allocation $(\rho_\gamma)$ exists.  For every $\gamma$, we have $\rho_\gamma\leq\alpha_\gamma^{-1}(p)$.
Each $\alpha_\gamma^{-1}(p)$ has rank one, so any non-zero
$\rho_\gamma$ must equal this orbit projection.  Since the orbit
projections are distinct and no two are orthogonal, the required
orthogonality of the $\rho_\gamma$'s allows at most one non-zero
term.  Such a term has rank one, and cannot sum to $I$.  This is a
contradiction.  Thus ergodicity and $p\neq0$ do not imply the
existence of a direct allocation, although
Corollary~\ref{cor:allocation-existence} gives a partial-isometric
allocation.
\end{example}

In the following proposition, we show that, in the classical setting, Theorem~\ref{thm:order-free-allocation} is the Kac formula for allocations of Meyerovitch and Weiss \cite[Proposition~1.3]{mw24}. Let $T\colon\Gamma\curvearrowright(X,\mathcal B,\mu)$ be a measure-preserving left action on a probability space, so that $T_{\gamma\eta}=T_\gamma\circ T_\eta$ for $\gamma,\eta\in\Gamma$. Let $\M=L^\infty(X,\mu)$, $\varphi(f)=\int_Xf~d\mu$ and $\alpha_\gamma(f)=f\circ T_\gamma^{-1}$. For a set $A\in\mathcal B$ with $\mu(A)>0$, let $p=\one_A$.

For an $A$-allocation $\kappa$, let $\rho_\gamma:=\one_{\kappa^{-1}(\gamma)}$ and $p_\gamma:=u_\gamma:=\alpha_\gamma(\rho_\gamma)$. For $z\in X$, let $B_\kappa(z):=\{\gamma\in\Gamma:\kappa(T_\gamma^{-1}(z))=\gamma\}$.

\begin{proposition}
\label{prop:MW-comparison}
For the action $\alpha_\gamma(f)=f\circ T_\gamma^{-1}$ and projection $p=\one_A$, the following hold.
\begin{enumerate}
\item Every partial-isometric $p$-allocation is direct. A family $(\rho_\gamma)_{\gamma\in\Gamma}$ is a direct $p$-allocation if and only if $\rho_\gamma=\one_{\kappa^{-1}(\gamma)}$ for every $\gamma\in\Gamma$, for some $A$-allocation $\kappa\colon X\to\Gamma$.
\item For every $A$-allocation $\kappa$, the set $B_\kappa(z)$ is finite for almost every $z\in A$. Also, $$\int_A|B_\kappa(z)|~d\mu(z)=1.$$ For $f\in L^\infty(X,\mu)$, we have $\int_AL_{\rho,p,u}(f)~d\mu=\int_Xf~d\mu$, and for almost every $z\in A$, we have
\begin{equation}
\label{eq:MW-fibre}
L_{\rho,p,u}(f)(z)=\sum_{\gamma\in B_\kappa(z)}f\bigl(T_\gamma^{-1}(z)\bigr).
\end{equation}
\item In the situation of \textup{(ii)}, for every measurable function $f\colon X\to[0,\infty]$, we have $$\int_A\sum_{\gamma\in B_\kappa(z)}f(T_\gamma^{-1}(z))~d\mu(z)=\int_Xf~d\mu.$$
\end{enumerate}
\end{proposition}

\begin{proof}
Theorem~\ref{thm:order-free-allocation} applies without ergodicity. For $E\in\mathcal B$ and $\gamma\in\Gamma$, we have $\alpha_\gamma(\one_E)=\one_E\circ T_\gamma^{-1}=\one_{T_\gamma(E)}$.

(i) By Remark~\ref{rem:direct-allocation}, every partial-isometric $p$-allocation for the abelian algebra $\M$ is direct. Let $(\rho_\gamma)_{\gamma\in\Gamma}$ be a direct $p$-allocation, and write $\rho_\gamma=\one_{E_\gamma}$. Since the $\rho_\gamma$ are pairwise orthogonal with sum $1$, the sets $E_\gamma$ are pairwise disjoint and cover $X$ modulo null sets; modifying them on null sets, we may assume that they partition $X$. Define $\kappa:=\gamma$ on $E_\gamma$. The inequality $\alpha_\gamma(\rho_\gamma)\leq p$ means that $\mu(T_\gamma(E_\gamma)\setminus A)=0$, and since $T_\gamma$ preserves $\mu$, this means that $T_\gamma(y)\in A$ for almost every $y\in E_\gamma$. As $\Gamma$ is countable, $T_{\kappa(y)}(y)\in A$ for almost every $y\in X$, so $\kappa$ is an $A$-allocation with $\rho_\gamma=\one_{\kappa^{-1}(\gamma)}$. Conversely, if $\kappa$ is an $A$-allocation, then the sets $\kappa^{-1}(\gamma)$ partition $X$, and since $T_\gamma(y)\in A$ for almost every $y\in\kappa^{-1}(\gamma)$, we have $\alpha_\gamma(\one_{\kappa^{-1}(\gamma)})=\one_{T_\gamma(\kappa^{-1}(\gamma))}\leq\one_A$.

(ii) For $\gamma\in\Gamma$ and $z\in X$, we have $z\in T_\gamma(\kappa^{-1}(\gamma))$ if and only if $\gamma\in B_\kappa(z)$. Since $\M$ is abelian and $u_\gamma=\alpha_\gamma(\rho_\gamma)$, it follows that
\[
u_\gamma^*\alpha_\gamma(f)u_\gamma=\alpha_\gamma(\rho_\gamma f)=\bigl(\one_{\kappa^{-1}(\gamma)}f\bigr)\circ T_\gamma^{-1}
\]
takes the value $f(T_\gamma^{-1}(z))$ at $z$ if $\gamma\in B_\kappa(z)$, and the value $0$ otherwise. Fix an enumeration of $\Gamma$. For $f=1$, the partial sums of $\sum_\gamma u_\gamma^*u_\gamma$ along this enumeration increase pointwise to the function $z\mapsto|B_\kappa(z)|$, and by Theorem~\ref{thm:order-free-allocation} they converge in $L^1(X,\mu)$ to $L_{\rho,p,u}(1)\in L^1(p\M p,\varphi)$, so a subsequence converges almost everywhere to $L_{\rho,p,u}(1)$. Hence $|B_\kappa(\cdot)|=L_{\rho,p,u}(1)$ almost everywhere, so $|B_\kappa(z)|$ is finite for almost every $z$ and vanishes for almost every $z\notin A$, and $\int_A|B_\kappa(z)|~d\mu(z)=\varphi(L_{\rho,p,u}(1))=\varphi(1)=1$.
Now let $f\in L^\infty(X,\mu)$. For almost every $z$, only the finitely many terms with $\gamma\in B_\kappa(z)$ are non-zero at $z$, so the partial sums along the enumeration converge almost everywhere to the right-hand side of \eqref{eq:MW-fibre}. By Theorem~\ref{thm:order-free-allocation}, they also converge in $L^1(X,\mu)$ to $L_{\rho,p,u}(f)$, so a subsequence converges almost everywhere to $L_{\rho,p,u}(f)$, which proves \eqref{eq:MW-fibre}. Moreover, $\int_AL_{\rho,p,u}(f)~d\mu=\varphi(L_{\rho,p,u}(f))=\varphi(f)=\int_Xf~d\mu$.

(iii) For $n\in\N$, the function $f\wedge n$ is bounded, and (ii) gives the asserted identity for $f\wedge n$. As $n\to\infty$, the functions $z\mapsto\sum_{\gamma\in B_\kappa(z)}(f\wedge n)(T_\gamma^{-1}(z))$ increase pointwise to $z\mapsto\sum_{\gamma\in B_\kappa(z)}f(T_\gamma^{-1}(z))$, and $f\wedge n$ increases to $f$. The monotone convergence theorem, applied on both sides, gives (iii).
\end{proof}

In the notation of \cite[Proposition~1.3]{mw24}, $f_\kappa(z)=\sum_{y\in X:~T_{\kappa(y)}(y)=z}f(y)$. The map $y\mapsto\kappa(y)$ is a bijection from $\{y\in X:T_{\kappa(y)}(y)=z\}$ onto $B_\kappa(z)$, with inverse $\gamma\mapsto T_\gamma^{-1}(z)$. Hence the right-hand side of \eqref{eq:MW-fibre} equals $f_\kappa(z)$, the identity in (iii) is \cite[(5)]{mw24}, and the identity $\int_A|B_\kappa(z)|~d\mu(z)=1$ in (ii) is \cite[(6)]{mw24}.

\subsection{\texorpdfstring{A center-valued coupling lemma}{A center-valued coupling lemma}}
\label{subsec:coupling}

The construction in Section~\ref{subsec:universal-towers} is modelled on the following classical picture. In the situation of Proposition~\ref{prop:MW-comparison}\textup{(ii)}, let $E_\gamma:=\kappa^{-1}(\gamma)$ for $\gamma\in\Gamma$ and $A_F:=\{z\in A:B_\kappa(z)=F\}$ for every finite set $F\subseteq\Gamma$. Since $B_\kappa(z)$ is finite for almost every $z\in A$, the sets $A_F$ partition $A$ modulo a null set. For $z\in X$, we have $z\in T_\gamma(E_\gamma)$ if and only if $\gamma\in B_\kappa(z)$, and since $T_\gamma(E_\gamma)\subseteq A$ modulo null sets,
\[
T_\gamma(E_\gamma)=\{z\in A:\gamma\in B_\kappa(z)\}=\biguplus_{F\ni\gamma}A_F~\text{modulo null sets}.
\]
Thus the target sets $T_\gamma(E_\gamma)$ of an allocation are unions of cells of the single partition $(A_F)_F$ of $A$. In the non-commutative setting, the analogous identities $p_\gamma=\sum_{F\ni\gamma}q_F$ for a single family $(q_F)_F$ of pairwise orthogonal projections force the projections $p_\gamma$ to commute, because every such sum belongs to the abelian von Neumann algebra generated by $(q_F)_F$. The projections $p_j$ of Lemma~\ref{lem:greedy-quantum-allocation} are not assumed to commute. We therefore first solve the corresponding problem for their center-valued traces, which lie in the abelian algebra $Z(\M)$ (Lemma~\ref{lem:finite-shape-threshold-coupling}), and then realize the resulting central elements by projections and partial isometries (Proposition~\ref{prop:universal-kac-towers1}).

No action is involved in this subsection: $(\M,\tau)$ is a finite von Neumann algebra with a faithful normal tracial state $\tau$, $\Theta\colon\M\to Z(\M)$ is the normalized center-valued trace of $\M$, and $z_{\mathrm{II}},z_1,z_2,\ldots$ are the central projections of its type decomposition (Section~\ref{subsec:center-valued-dimension-tools}). We keep the enumeration $(g_j)_{j\in J}$ of $\Gamma$, and we write $\mathcal F$ for the set of finite subsets of $\Gamma$ containing $e$.

For projections $p,p_j\in\M$, $j\in J$, write $a:=\Theta(p)$ and $\delta_j:=\Theta(p_j)$.

\begin{lemma}
\label{lem:finite-shape-threshold-coupling}
Let $p\in\M$ be a non-zero projection, and let $(p_j)_{j\in J}$ be projections in $\M$ with $p_1=p$, $p_j\leq p$ for every $j\in J$, and $\sum_{j\in J}\tau(\Theta(p_j))<\infty$. Then there is a family $(c_F)_{F\in\mathcal F}$ in $Z(\M)_+$ such that, for every $j\in J$, we have $\sum_{F\in\mathcal F}c_F=a$ and $\sum_{F\in\mathcal F,~F\ni g_j}c_F=\delta_j$, where both sums converge strongly. Moreover, for all $m\geq1$ and $F\in\mathcal F$, the element $mz_mc_F$ is integer-valued in $Z(z_m\M)$, in the sense of Lemma~\ref{lem:type-I-dimension-range}. In particular, every enumeration of $(c_F)_{F\in\mathcal F}$, followed by zero terms if $\mathcal F$ is finite, satisfies the hypotheses of Corollary~\ref{cor:center-valued-dimension-partition} for the projection $p$.
\end{lemma}

\begin{proof}
The GNS representation of $Z(\M)$ associated with $\tau|_{Z(\M)}$ has a cyclic and separating unit vector, since $\tau$ is faithful. Applying \cite[Proposition~III.1.5.16]{MR2188261} in this representation, we identify $Z(\M)$ with $L^\infty(\Omega,\mu)$ for a probability space $(\Omega,\mu)$ such that $\tau(z)=\int_\Omega z~d\mu$ for $z\in Z(\M)$. Since $\Theta$ is positive and $p_j\leq p=p_1$, we have $0\leq\delta_j\leq a$ and $\delta_1=a$. As $J$ is countable, we may choose representatives of $a$ and of the $\delta_j$ which satisfy these relations everywhere on $\Omega$.

Let $S:=\{\omega\in\Omega:a(\omega)>0\}$. On $S$, define $r_j:=\delta_j/a$. On $\Omega\setminus S$, let $r_1:=1$ and $r_j:=0$ for $j\neq1$. Then $0\leq r_j\leq1$, $r_1=1$, and $ar_j=\delta_j$ on $\Omega$, because $0\leq\delta_j\leq a$ forces $\delta_j=0$ on $\Omega\setminus S$. By Tonelli's theorem, $\int_\Omega\sum_{j\in J}\delta_j~d\mu=\sum_{j\in J}\tau(\delta_j)<\infty$, so $\sum_j\delta_j(\omega)<\infty$, and hence $\sum_jr_j(\omega)<\infty$, for almost every $\omega\in S$. On $\Omega\setminus S$, we have $\sum_jr_j=1$. Redefining $r_1:=1$ and $r_j:=0$ for $j\neq1$ on the exceptional null set, we may assume that $\sum_jr_j(\omega)<\infty$ for every $\omega\in\Omega$, while $ar_j=\delta_j$ still holds almost everywhere.

For $\omega\in\Omega$ and $0<t\leq1$, define $F(\omega,t):=\{g_j:j\in J,~r_j(\omega)\geq t\}$. This set contains $e=g_1$, since $r_1=1$, and it is finite: each of its elements contributes at least $t$ to $\sum_jr_j(\omega)$, so it has at most $t^{-1}\sum_jr_j(\omega)$ elements. Thus $F(\omega,t)\in\mathcal F$. For $F\in\mathcal F$, let
\[
G_F:=\{(\omega,t)\in\Omega\times(0,1]:F(\omega,t)=F\}
=\bigcap_{g_j\in F}\{(\omega,t):t\leq r_j(\omega)\}\cap\bigcap_{g_j\notin F}\{(\omega,t):t>r_j(\omega)\}.
\]
This is a countable intersection of measurable subsets of $\Omega\times(0,1]$, and hence measurable. By Tonelli's theorem, $\lambda_F(\omega):=\int_0^1\one_{G_F}(\omega,t)~dt$ defines a measurable function on $\Omega$ with $0\leq\lambda_F\leq1$. Let $c_F:=a\lambda_F\in Z(\M)_+$.

The set $\mathcal F$ is countable, and for every $(\omega,t)$ exactly one $F\in\mathcal F$ satisfies $F(\omega,t)=F$. Hence $\sum_{F\in\mathcal F}\one_{G_F}=1$ on $\Omega\times(0,1]$, and integrating in $t$ (Tonelli) gives $\sum_{F\in\mathcal F}\lambda_F(\omega)=1$ for every $\omega\in\Omega$. Moreover, for $j\in J$ and $\omega\in\Omega$,
\[
\sum_{F\in\mathcal F,~F\ni g_j}\lambda_F(\omega)
=\int_0^1\sum_{F\ni g_j}\one_{G_F}(\omega,t)~dt
=\int_0^1\one_{\{t\leq r_j(\omega)\}}~dt
=r_j(\omega).
\]
Multiplying by $a$, we obtain $\sum_Fc_F=a$ everywhere and $\sum_{F\ni g_j}c_F=ar_j=\delta_j$ almost everywhere. In both series the terms are non-negative and bounded by $a$, so their finite partial sums form increasing nets which converge pointwise, and by dominated convergence also in the weak$^*$ topology of $L^\infty(\Omega,\mu)$, to the sums computed above. By \cite[Section~1.1, p.~1]{hiai2021}, each bounded increasing net converges strongly to its supremum, so both series converge strongly.

It remains to check integrality. Fix $m\geq1$ with $z_m\neq0$, and write $z_m=\one_{\Omega_m}$. For every projection $q\in\M$, the element $z_mq$ is a projection in the homogeneous type~$\mathrm I_m$ algebra $z_m\M$, whose normalized center-valued trace is the restriction of $\Theta$. Lemma~\ref{lem:type-I-dimension-range}, applied in $z_m\M$ with $r=z_mq$ and $b=\Theta(z_mq)$, shows that $m\Theta(z_mq)=mz_m\Theta(q)$ is integer-valued. Hence, after discarding a null set, the functions $k:=ma$ and $\ell_j:=m\delta_j$ take values in $\{0,1,\ldots,m\}$ on $\Omega_m$, and $\ell_j\leq k$ there. Fix $\omega\in\Omega_m$. If $k(\omega)=0$, then $a(\omega)=0$ and $c_F(\omega)=0$ for every $F$. If $k(\omega)>0$, then $r_j(\omega)=\ell_j(\omega)/k(\omega)$ belongs to $\{0,1/k(\omega),\ldots,1\}$ for every $j$. Consequently, for $h\in\{1,\ldots,k(\omega)\}$ and $t\in((h-1)/k(\omega),h/k(\omega)]$, we have $r_j(\omega)\geq t$ if and only if $r_j(\omega)\geq h/k(\omega)$; thus $t\mapsto F(\omega,t)$ is constant on each such interval, and $\lambda_F(\omega)$ is an integer multiple of $1/k(\omega)$. Therefore $mc_F(\omega)=k(\omega)\lambda_F(\omega)$ is an integer in $\{0,\ldots,m\}$. The projections $e_i:=\one_{\Omega_m\cap\{mc_F=i\}}$, $i=0,\ldots,m$, are pairwise orthogonal central projections with $\sum_ie_i=z_m$ and $mz_mc_F=\sum_iie_i$, as required.
\end{proof}

\subsection{Universal partial-isometric Kac towers}
\label{subsec:universal-towers}

In this subsection, we generalise the universal-shape theorem of Meyerovitch and Weiss \cite[Theorem~1.5]{mw24} to actions on finite von Neumann algebras. Ergodicity is not assumed. Throughout, $(\M,\tau)$, $\Theta$ and $\mathcal F$ are as in Section~\ref{subsec:coupling}, $\alpha\colon\Gamma\to\Aut(\M)$ is a $\tau$-preserving action, and $p\in\M$ is a non-zero sweep-out projection, that is, $\bigvee_{\gamma\in\Gamma}\alpha_\gamma(p)=1$. If $\alpha$ is ergodic and preserves a faithful normal state $\varphi$, then these assumptions hold with $\tau=\varphi$ for every non-zero projection $p$ (Corollary~\ref{cor:universal-towers-ergodic}).

\begin{definition}
\label{def:kac-tower}
Let $(B_n)_{n\geq1}$ be a sequence of finite subsets of $\Gamma$, each containing $e$. A partial-isometric Kac tower for $(\alpha,p)$ with shapes $(B_n)_{n\geq1}$ consists of projections $q_n\in\Proj(p\M p)$, projections $\rho_{n,\gamma}\in\Proj(\M)$ and partial isometries $v_{n,\gamma}\in\M$, for $n\geq1$ and $\gamma\in B_n$, with the following properties.
\begin{enumerate}
\item The projections $q_n$, $n\geq1$, are pairwise orthogonal, and $\sum_{n\geq1}q_n=p$.
\item The projections $\rho_{n,\gamma}$, $n\geq1$, $\gamma\in B_n$, are pairwise orthogonal, and $\sum_{n\geq1}\sum_{\gamma\in B_n}\rho_{n,\gamma}=1$.
\item For all $n\geq1$ and $\gamma\in B_n$, we have $v_{n,\gamma}^*v_{n,\gamma}=q_n$ and $v_{n,\gamma}v_{n,\gamma}^*=\alpha_\gamma(\rho_{n,\gamma})$.
\item For every $n\geq1$, we have $\rho_{n,e}=v_{n,e}=q_n$.
\end{enumerate}
A pair $(n,\gamma)$ with $n\geq1$ and $\gamma\in B_n$ is called a branch. When $\M=L^\infty(X,\mu)$ and $p=\one_A$ as in Proposition~\ref{prop:commutative-universal-towers}, we have $q_n=\one_{A_n}$ for a measurable set $A_n\subseteq A$ and $\rho_{n,\gamma}=\one_{T_\gamma^{-1}(A_n)}$ for $\gamma\in B_n$, so that the $\rho_{n,\gamma}$, $\gamma\in B_n$, correspond to the levels of a tower over $A_n$ with shape $B_n$.
\end{definition}

\begin{proposition}
\label{prop:universal-kac-towers1}
There is a sequence $(B_n)_{n\geq1}$ of finite subsets of $\Gamma$, each containing $e$, which depends only on $\Gamma$ and has the following property: for every finite von Neumann algebra $(\M,\tau)$ with a faithful normal tracial state, every $\tau$-preserving action $\alpha\colon\Gamma\to\Aut(\M)$ and every non-zero sweep-out projection $p\in\M$, there is a partial-isometric Kac tower for $(\alpha,p)$ with shapes $(B_n)_{n\geq1}$.
\end{proposition}

\begin{proof}
We first choose the shapes. Since $\Gamma$ is countable, so is $\mathcal F$, and $\{e\}\in\mathcal F$. Choose a sequence $(B_n)_{n\geq1}$ in which every member of $\mathcal F$ occurs; if $\Gamma$ is finite, then $\mathcal F$ is finite and repetitions are necessary. For each $F\in\mathcal F$, choose an index $\iota(F)$ with $B_{\iota(F)}=F$. Then $\iota\colon\mathcal F\to\N$ is injective, and these choices depend only on $\Gamma$.

Now let $(\M,\tau)$, $\alpha$ and $p$ be as in the statement. By Lemma~\ref{lem:greedy-quantum-allocation}, there are pairwise orthogonal projections $(e_j)_{j\in J}$ with $\sum_je_j=1$, and projections $d_j\leq\alpha_{g_j}^{-1}(p)$ with $e_j\sim d_j$ and $e_1=d_1=p$. Let $p_j:=\alpha_{g_j}(d_j)$ and $\delta_j:=\Theta(p_j)$; thus $p_j\leq p$ and $p_1=p$. For every $j\in J$, the equivariance $\Theta\circ\alpha_\gamma=\alpha_\gamma\circ\Theta$ (Section~\ref{subsec:center-valued-dimension-tools}) and the equivalence $d_j\sim e_j$ give
\begin{equation}
\label{eq:tower-twisted-dimensions}
\alpha_{g_j}^{-1}(\delta_j)=\Theta\bigl(\alpha_{g_j}^{-1}(p_j)\bigr)=\Theta(d_j)=\Theta(e_j).
\end{equation}
By \cite[Theorem~III.2.5.7\textup{(iv)}]{MR2188261}, $\tau\circ\Theta=\tau$. Traciality, invariance and normality of $\tau$ give
\[
\sum_{j\in J}\tau(\Theta(p_j))=\sum_{j\in J}\tau(p_j)=\sum_{j\in J}\tau(d_j)=\sum_{j\in J}\tau(e_j)=1.
\]
Thus Lemma~\ref{lem:finite-shape-threshold-coupling} applies to $p$ and $(p_j)_{j\in J}$, and it gives central elements $(c_F)_{F\in\mathcal F}$ such that, for every $j\in J$,
\begin{equation}
\label{eq:tower-central-marginals}
\sum_{F\in\mathcal F}c_F=\Theta(p)~\text{and}~\sum_{F\in\mathcal F,~F\ni g_j}c_F=\delta_j.
\end{equation}
These elements satisfy the integrality hypothesis of Corollary~\ref{cor:center-valued-dimension-partition}. Applying that corollary to $p$, we obtain pairwise orthogonal projections $q_F\leq p$, $F\in\mathcal F$, with $\Theta(q_F)=c_F$ and $\sum_{F\in\mathcal F}q_F=p$.

We next refine the projections $e_j$. Fix $j\in J$ with $j\neq1$, and let $\mathcal F_j:=\{F\in\mathcal F:g_j\in F\}$. Since $\alpha_{g_j}^{-1}$ is normal, \eqref{eq:tower-central-marginals} and \eqref{eq:tower-twisted-dimensions} give $\sum_{F\in\mathcal F_j}\alpha_{g_j}^{-1}(c_F)=\alpha_{g_j}^{-1}(\delta_j)=\Theta(e_j)$.
Automorphisms preserve abelian projections and Murray--von Neumann equivalence, hence preserve homogeneity of type~$\mathrm I_m$. The canonical central type decomposition \cite[Sections~III.1.4.7 and~III.1.5.12--III.1.5.13]{MR2188261} therefore gives, for every $m\geq1$, $\alpha_{g_j}^{-1}(z_m)=z_m$. Consequently, if $mz_mc_F=\sum_{i=0}^mi~z_{m,i}$ with pairwise orthogonal central projections $z_{m,i}$ summing to $z_m$, then $mz_m\alpha_{g_j}^{-1}(c_F)=\sum_{i=0}^mi~\alpha_{g_j}^{-1}(z_{m,i})$, where the $\alpha_{g_j}^{-1}(z_{m,i})$ are pairwise orthogonal central projections summing to $z_m$. Thus the family $(\alpha_{g_j}^{-1}(c_F))_{F\in\mathcal F_j}$ satisfies the hypotheses of Corollary~\ref{cor:center-valued-dimension-partition} for the projection $e_j$, and we obtain pairwise orthogonal projections $\rho_{F,g_j}\leq e_j$, $F\in\mathcal F_j$, with $\Theta(\rho_{F,g_j})=\alpha_{g_j}^{-1}(c_F)$ and $\sum_{F\in\mathcal F_j}\rho_{F,g_j}=e_j$.
For the identity label, we let $\rho_{F,e}:=q_F$ for $F\in\mathcal F$; since every $F\in\mathcal F$ contains $e$ and $\sum_Fq_F=p=e_1$, these projections partition $e_1$. Projections with different labels lie under different $e_j$'s and are therefore orthogonal. Hence the projections $\rho_{F,\gamma}$, $F\in\mathcal F$, $\gamma\in F$, are pairwise orthogonal, and their sum is $\sum_{j\in J}e_j=1$.

Next, for $F\in\mathcal F$ and $g_j\in F$ with $j\neq1$, equivariance gives $\Theta\bigl(\alpha_{g_j}(\rho_{F,g_j})\bigr)=\alpha_{g_j}\bigl(\Theta(\rho_{F,g_j})\bigr)=c_F=\Theta(q_F)$, so the comparison theorem \cite[Theorem~8.4.3(vi), p.~532]{MR1468230} gives a partial isometry $v_{F,g_j}\in\M$ with $v_{F,g_j}^*v_{F,g_j}=q_F$ and $v_{F,g_j}v_{F,g_j}^*=\alpha_{g_j}(\rho_{F,g_j})$. For the identity label, we let $v_{F,e}:=q_F$, which satisfies the same relations because $\rho_{F,e}=q_F$ and $\alpha_e=\mathrm{id}$.

Finally, we index these objects by the sequence $(B_n)_{n\geq1}$. For $F\in\mathcal F$ and $\gamma\in F$, let $q_{\iota(F)}:=q_F$, $\rho_{\iota(F),\gamma}:=\rho_{F,\gamma}$ and $v_{\iota(F),\gamma}:=v_{F,\gamma}$. For $n\notin\iota(\mathcal F)$, let $q_n:=0$ and $\rho_{n,\gamma}:=v_{n,\gamma}:=0$ for $\gamma\in B_n$. Since $\iota$ is injective, properties (i) and (ii) of Definition~\ref{def:kac-tower} follow from $\sum_Fq_F=p$ and from the previous two paragraphs, while (iii) and (iv) hold by construction, and trivially for $n\notin\iota(\mathcal F)$.
\end{proof}

\begin{remark}
\label{rem:partial-isometries-essential}
The partial isometries in Definition~\ref{def:kac-tower} cannot be replaced by equalities of projections. Indeed, suppose that a partial-isometric Kac tower satisfies $\alpha_\gamma(\rho_{n,\gamma})=q_n$ for every branch $(n,\gamma)$, and let $\eta_\gamma:=\sum_{n\geq1,~\gamma\in B_n}\rho_{n,\gamma}$ for $\gamma\in\Gamma$. By Definition~\ref{def:kac-tower}\textup{(ii)}, the $\eta_\gamma$ are pairwise orthogonal projections with $\sum_\gamma\eta_\gamma=1$, and by normality of $\alpha_\gamma$ and Definition~\ref{def:kac-tower}\textup{(i)}, $\alpha_\gamma(\eta_\gamma)=\sum_{n\geq1,~\gamma\in B_n}q_n\leq p$. Thus $(\eta_\gamma)_{\gamma\in\Gamma}$ is a direct $p$-allocation. For the action and the projection of Example~\ref{ex:direct-allocation-fails}, no direct $p$-allocation exists. Hence no partial-isometric Kac tower for them satisfies $\alpha_\gamma(\rho_{n,\gamma})=q_n$ for all branches, although Proposition~\ref{prop:universal-kac-towers1} provides a partial-isometric Kac tower.
\end{remark}

For the rest of this subsection, we fix a sequence $(B_n)_{n\geq1}$ of finite subsets of $\Gamma$ containing $e$ and a partial-isometric Kac tower $(q_n,\rho_{n,\gamma},v_{n,\gamma})$ for $(\alpha,p)$ with shapes $(B_n)_{n\geq1}$, for instance the one given by Proposition~\ref{prop:universal-kac-towers1}. We first use the partial isometries $v_{n,\gamma}$ to identify each source corner with its target corner.

For a branch $(n,\gamma)$, write $q=q_n$, $\rho=\rho_{n,\gamma}$ and $v=v_{n,\gamma}$, and define $\beta_{n,\gamma}\colon\rho\M\rho\longrightarrow q\M q$ by $\beta_{n,\gamma}(a)=v^*\alpha_\gamma(a)v$.

\begin{proposition}
\label{prop:universal-kac-towers2}
For every branch $(n,\gamma)$, the map $\beta_{n,\gamma}$ is a normal $*$-isomorphism with $\beta_{n,\gamma}(\rho)=q$ and $\tau\circ\beta_{n,\gamma}=\tau$ on $\rho\M\rho$. It extends uniquely to an isometry of $L^1(\rho\M\rho,\tau|_{\rho\M\rho})$ onto $L^1(q\M q,\tau|_{q\M q})$, again denoted by $\beta_{n,\gamma}$, which satisfies $\tau\circ\beta_{n,\gamma}=\tau$. Moreover, $\alpha_\gamma$ extends uniquely to an isometry of $L^1(\M,\tau)$, again denoted by $\alpha_\gamma$. For $x\in L^1(\M,\tau)$, we have
\begin{equation}
\label{eq:branch-map}
v^*\alpha_\gamma(x)v=\beta_{n,\gamma}(\rho x\rho).
\end{equation}
\end{proposition}

\begin{proof}
Since $vv^*=\alpha_\gamma(\rho)$, we have $v=\alpha_\gamma(\rho)v$. For $a,b\in\rho\M\rho$,
\[
\beta_{n,\gamma}(a)\beta_{n,\gamma}(b)=v^*\alpha_\gamma(a)vv^*\alpha_\gamma(b)v=v^*\alpha_\gamma(a\rho b)v=\beta_{n,\gamma}(ab),
\]
and clearly $\beta_{n,\gamma}(a^*)=\beta_{n,\gamma}(a)^*$ and $\beta_{n,\gamma}(\rho)=v^*vv^*v=q$. For $b\in q\M q$, we have $vbv^*=\alpha_\gamma(\rho)vbv^*\alpha_\gamma(\rho)$, so $b\mapsto\alpha_\gamma^{-1}(vbv^*)$ maps $q\M q$ into $\rho\M\rho$. It is inverse to $\beta_{n,\gamma}$, because, for $a\in\rho\M\rho$ and $b\in q\M q$, we have $\alpha_\gamma^{-1}(vv^*\alpha_\gamma(a)vv^*)=\rho a\rho=a$ and $v^*vbv^*v=qbq=b$.
Both maps are compositions of automorphisms and of maps of the form $y\mapsto w^*yw$, and hence normal. By traciality and invariance of $\tau$, for $a\in\rho\M\rho$, we have $\tau\bigl(v^*\alpha_\gamma(a)v\bigr)=\tau\bigl(\alpha_\gamma(a)vv^*\bigr)=\tau\bigl(\alpha_\gamma(a\rho)\bigr)=\tau(a)$.
A $*$-isomorphism preserves positive square roots, so $|\beta_{n,\gamma}(a)|=\beta_{n,\gamma}(|a|)$, and therefore $\|\beta_{n,\gamma}(a)\|_1=\tau(\beta_{n,\gamma}(|a|))=\tau(|a|)=\|a\|_1$.
For $e=\rho,q$, we identify $L^1(e\M e,\tau|_{e\M e})$ with $eL^1(\M,\tau)e$: the inclusion of the bounded corner preserves the $L^1$-norm, and its extension has this range by density and boundedness of compression \cite[Theorem~4.41 and Lemma~4.35(2)]{hiai2021}. By \cite[Theorem~4.41]{hiai2021}, $\rho\M\rho$ is dense in $L^1(\rho\M\rho,\tau|_{\rho\M\rho})$, so $\beta_{n,\gamma}$ extends uniquely to an isometry of this space into $L^1(q\M q,\tau|_{q\M q})$. Its range is closed and contains the dense subspace $\beta_{n,\gamma}(\rho\M\rho)=q\M q$, so the extension is onto. Since $|\tau(y)|\leq\|y\|_1$ by \cite[Proposition~4.37]{hiai2021}, the identity $\tau\circ\beta_{n,\gamma}=\tau$ extends by continuity. The same argument, applied to the $\tau$-preserving automorphism $\alpha_\gamma$ of $\M$, gives the isometric extension of $\alpha_\gamma$ to $L^1(\M,\tau)$. Finally, for $x\in\M$,
\[
v^*\alpha_\gamma(x)v=v^*\alpha_\gamma(\rho)\alpha_\gamma(x)\alpha_\gamma(\rho)v=v^*\alpha_\gamma(\rho x\rho)v=\beta_{n,\gamma}(\rho x\rho).
\]
By \cite[Lemma~4.35(2)]{hiai2021}, both sides of \eqref{eq:branch-map} are continuous in $x\in L^1(\M,\tau)$, and $\M$ is dense in $L^1(\M,\tau)$ by \cite[Theorem~4.41]{hiai2021}, so \eqref{eq:branch-map} holds for every $x\in L^1(\M,\tau)$.
\end{proof}

For $x\in L^1(\M,\tau)$, we consider the series
\[
L_p(x):=\sum_{n\geq1}\sum_{\gamma\in B_n}v_{n,\gamma}^*\alpha_\gamma(x)v_{n,\gamma}=\sum_{n\geq1}\sum_{\gamma\in B_n}\beta_{n,\gamma}(\rho_{n,\gamma}x\rho_{n,\gamma}),
\]
whose two expressions agree term by term by \eqref{eq:branch-map}.

\begin{proposition}
\label{prop:universal-kac-towers3}
For every $x\in L^1(\M,\tau)$, we have $\sum_{n\geq1}\sum_{\gamma\in B_n}\|v_{n,\gamma}^*\alpha_\gamma(x)v_{n,\gamma}\|_1\leq\|x\|_1$. In particular, the series $L_p(x)$ converges absolutely in $L^1(p\M p,\tau|_{p\M p})$. Furthermore, $\tau(L_p(x))=\tau(x)$. If $x\geq0$, then every term of the series is positive, and $L_p(x)\in L^1(p\M p,\tau|_{p\M p})_+$.
\end{proposition}

\begin{proof}
Let $K$ be a finite set of branches. Let $\rho_K:=\sum_{(n,\gamma)\in K}\rho_{n,\gamma}$ and, for $y\in L^1(\M,\tau)$, let $\mathbb E_K(y):=\sum_{(n,\gamma)\in K}\rho_{n,\gamma}y\rho_{n,\gamma}$. We first show that $\|\mathbb E_K(y)\|_1\leq\|y\|_1$. If $\rho_K=0$, then $\mathbb E_K=0$. Otherwise, on $\rho_K\M\rho_K$, orthogonality of the $\rho_{n,\gamma}$ shows that $\mathbb E_K$ fixes $\bigoplus_{(n,\gamma)\in K}\rho_{n,\gamma}\M\rho_{n,\gamma}$, is bimodular over this algebra, and preserves the normalized trace $\tau(\rho_K)^{-1}\tau|_{\rho_K\M\rho_K}$. Thus \cite[Theorem~5.4\textup{(1)--(2)}]{hiai2021} identifies it as the corresponding conditional expectation and gives its contractive extension to $L^1$. By \cite[Theorem~4.41 and Lemma~4.35(2)]{hiai2021}, this extension has the same compression formula. Since $\mathbb E_K(y)=\mathbb E_K(\rho_Ky\rho_K)$ and normalization multiplies both $L^1$-norms by the same constant, we obtain $\|\mathbb E_K(y)\|_1\leq\|\rho_Ky\rho_K\|_1\leq\|y\|_1$.

Next, let $x\in L^1(\M,\tau)$ and $w_{n,\gamma}:=\rho_{n,\gamma}x\rho_{n,\gamma}$. We show that $\sum_{(n,\gamma)\in K}\|w_{n,\gamma}\|_1=\|\mathbb E_K(x)\|_1$. The inequality $\geq$ is the triangle inequality. For the converse inequality, write $w_{n,\gamma}=s_{n,\gamma}|w_{n,\gamma}|$ by polar decomposition \cite[Section~4.1, p.~39]{hiai2021}. Since $(1-\rho_{n,\gamma})w_{n,\gamma}=0=w_{n,\gamma}(1-\rho_{n,\gamma})$, the initial and final projections of $s_{n,\gamma}$ lie under $\rho_{n,\gamma}$. Hence $s:=\sum_{(n,\gamma)\in K}s_{n,\gamma}^*$ satisfies $s^*s=\sum_{(n,\gamma)\in K}s_{n,\gamma}s_{n,\gamma}^*\leq\rho_K$, so $\|s\|\leq1$, and $s_{n,\gamma}^*w_{n',\gamma'}=0$ for $(n,\gamma)\neq(n',\gamma')$. Therefore \cite[Lemma~4.46]{hiai2021} gives
\[
\|\mathbb E_K(x)\|_1\geq\bigl|\tau\bigl(s~\mathbb E_K(x)\bigr)\bigr|=\sum_{(n,\gamma)\in K}\tau\bigl(s_{n,\gamma}^*w_{n,\gamma}\bigr)=\sum_{(n,\gamma)\in K}\|w_{n,\gamma}\|_1.
\]
Combining these two facts with \eqref{eq:branch-map} and the isometry property in Proposition~\ref{prop:universal-kac-towers2}, we obtain
\[
\sum_{(n,\gamma)\in K}\bigl\|v_{n,\gamma}^*\alpha_\gamma(x)v_{n,\gamma}\bigr\|_1=\sum_{(n,\gamma)\in K}\|w_{n,\gamma}\|_1=\|\mathbb E_K(x)\|_1\leq\|x\|_1.
\]
Taking the supremum over all finite sets $K$ proves the asserted inequality. As in the proof of Theorem~\ref{thm:order-free-allocation}, the finite partial sums $S_K:=\sum_{(n,\gamma)\in K}v_{n,\gamma}^*\alpha_\gamma(x)v_{n,\gamma}$ converge in $L^1(p\M p,\tau|_{p\M p})$ to $L_p(x)$.

By Proposition~\ref{prop:universal-kac-towers2} and \cite[Proposition~4.45]{hiai2021},
\[
\tau(S_K)=\sum_{(n,\gamma)\in K}\tau(\rho_{n,\gamma}x\rho_{n,\gamma})=\sum_{(n,\gamma)\in K}\tau(x\rho_{n,\gamma})=\tau(x\rho_K).
\]
The projections $\rho_K$ increase to $1$, and by \cite[Corollary~4.51]{hiai2021} the functional $y\mapsto\tau(xy)$ is normal on $\M$; hence $\tau(x\rho_K)\to\tau(x)$. Since $|\tau(y)|\leq\|y\|_1$ by \cite[Proposition~4.37]{hiai2021}, we obtain $\tau(L_p(x))=\lim_K\tau(S_K)=\tau(x)$.

Finally, let $x\geq0$. By \cite[Proposition~4.25]{hiai2021}, the operators $x\wedge k\in\M_+$ converge to $x$ in $L^1(\M,\tau)$. For every branch, $v_{n,\gamma}^*\alpha_\gamma(x\wedge k)v_{n,\gamma}\geq0$, and these operators converge in $L^1$ to $v_{n,\gamma}^*\alpha_\gamma(x)v_{n,\gamma}$ by \cite[Lemma~4.35(2)]{hiai2021} and Proposition~\ref{prop:universal-kac-towers2}. The positive cone is norm closed by \cite[Corollary~4.51]{hiai2021}, so every term of the series is positive. Hence every $S_K$ is positive, and its $L^1$-limit $L_p(x)$ is positive.
\end{proof}

For each branch $(n,\gamma)$, let $\Psi_{n,\gamma}\colon\M\to p\M p$ be the normal positive map defined, for $a\in\M$, by $\Psi_{n,\gamma}(a):=v_{n,\gamma}^*\alpha_\gamma(a)v_{n,\gamma}=\beta_{n,\gamma}(\rho_{n,\gamma}a\rho_{n,\gamma})$.
We use its extension $\widehat{\Psi_{n,\gamma}}$ from Section~\ref{subsec:extended-positive-part} to define the sum for extended positive inputs. For $x\in\widehat{\M}_+$, let $\widehat L_p(x):=\sup_{N\geq1}~\sum_{n=1}^N\sum_{\gamma\in B_n}\widehat{\Psi_{n,\gamma}}(x)\in\widehat{(p\M p)}_+$.
The finite sums increase with $N$, because each $\widehat{\Psi_{n,\gamma}}(x)$ belongs to $\widehat{(p\M p)}_+$, and so the supremum exists by \cite[Section~7.1, p.~95]{hiai2021}. For $m\in\widehat{\M}_+$, we let $\widehat\tau(m):=m(\tau)$, and for $m\in\widehat{(p\M p)}_+$, we let $\widehat\tau(m):=m(\tau|_{p\M p})$.

Let $R_p:=L_p(1)$ and $\Phi_p:=\sum_{n\geq1}nq_n$, where the latter is a positive spectral series. For $n\geq1$, let $h(n):=|B_n|$, and let $\tau_p(y):=\tau(y)/\tau(p)$ be the normalized corner trace.

\begin{theorem}
\label{thm:universal-kac-towers}
\leavevmode
\begin{enumerate}
\item For every $x\in\widehat{\M}_+$, we have $\widehat\tau(\widehat L_p(x))=\widehat\tau(x)$ in $[0,\infty]$. For every $x\in L^1(\M,\tau)_+$, we have $\widehat L_p(\widehat x)=\widehat{L_p(x)}$.
\item The operator $R_p=\sum_{n\geq1}|B_n|q_n$ belongs to $L^1(p\M p,\tau|_{p\M p})_+$. Moreover, $$\tau(R_p)=\sum_{n\geq1}|B_n|~\tau(q_n)=1.$$ Consequently, $\tau_p(R_p)=1/\tau(p)$.
\item The positive self-adjoint operator $\Phi_p$ is affiliated with $p\M p$, and for each $n\geq1$, we have $\one_{\{n\}}(\Phi_p)=q_n$. Moreover, $R_p=h(\Phi_p)$.
\end{enumerate}
\end{theorem}

\begin{proof}
(i) Let $x\in\widehat{\M}_+$. For a branch $(n,\gamma)$ and $y\in\M$, Proposition~\ref{prop:universal-kac-towers2} and traciality give $\tau(\Psi_{n,\gamma}(y))=\tau(\rho_{n,\gamma}y\rho_{n,\gamma})=\tau(y\rho_{n,\gamma})$, that is, $\tau|_{p\M p}\circ\Psi_{n,\gamma}=\tau(~\cdot~\rho_{n,\gamma})$. For $N\geq1$, let $\rho_{(N)}:=\sum_{n\leq N}\sum_{\gamma\in B_n}\rho_{n,\gamma}$ and $\omega_N:=\tau(~\cdot~\rho_{(N)})$. By the definition of $\widehat{\Psi_{n,\gamma}}$ and the additivity of $x$,
\[
\widehat\tau\Bigl(\sum_{n\leq N}\sum_{\gamma\in B_n}\widehat{\Psi_{n,\gamma}}(x)\Bigr)=\sum_{n\leq N}\sum_{\gamma\in B_n}x\bigl(\tau(~\cdot~\rho_{n,\gamma})\bigr)=x(\omega_N).
\]
The functional $\tau-\omega_N=\tau(~\cdot~(1-\rho_{(N)}))$ is positive, since $\tau(y(1-\rho_{(N)}))=\tau((1-\rho_{(N)})y(1-\rho_{(N)}))$, and $\|\tau-\omega_N\|=\tau(1-\rho_{(N)})\to0$ by normality and Definition~\ref{def:kac-tower}\textup{(ii)}. By additivity and lower semicontinuity of $x$ \cite[Definition~7.1]{hiai2021}, we have $x(\omega_N)\leq x(\tau)\leq\liminf_Nx(\omega_N)$, so $x(\tau)=\lim_Nx(\omega_N)$. Since $\widehat\tau$ commutes with suprema of increasing sequences by \cite[Proposition~7.5]{hiai2021}, we conclude that $\widehat\tau(\widehat L_p(x))=\sup_Nx(\omega_N)=x(\tau)=\widehat\tau(x)$.

Now let $x\in L^1(\M,\tau)_+$. By definition, $\widehat x=\sup_k(x\wedge k)$, and by \cite[Proposition~4.25]{hiai2021}, $x\wedge k\to x$ in $L^1(\M,\tau)$. For a branch $(n,\gamma)$, the operators $\Psi_{n,\gamma}(x\wedge k)$ increase with $k$ and converge in $L^1$ to $v_{n,\gamma}^*\alpha_\gamma(x)v_{n,\gamma}$, by \cite[Lemma~4.35(2)]{hiai2021} and Proposition~\ref{prop:universal-kac-towers2}. By the definition of $\widehat{\Psi_{n,\gamma}}$, evaluating at each positive normal functional shows that it commutes with suprema of increasing nets. Applying \cite[Proposition~4.23(2)]{hiai2021} as in the proof of Theorem~\ref{thm:order-free-allocation}, we obtain $\widehat{\Psi_{n,\gamma}}(\widehat x)=\sup_k\Psi_{n,\gamma}(x\wedge k)=\widehat{v_{n,\gamma}^*\alpha_\gamma(x)v_{n,\gamma}}$.
For $N\geq1$, let $S_N:=\sum_{n\leq N}\sum_{\gamma\in B_n}v_{n,\gamma}^*\alpha_\gamma(x)v_{n,\gamma}$. By \cite[Corollary~4.24]{hiai2021}, applied with each positive $L^1$-density, and \cite[Example~7.2(1) and Corollary~4.51]{hiai2021}, we have $\sum_{n\leq N}\sum_{\gamma\in B_n}\widehat{\Psi_{n,\gamma}}(\widehat x)=\widehat{S_N}$. By Proposition~\ref{prop:universal-kac-towers3}, $(S_N)_N$ is an increasing sequence in $L^1(p\M p,\tau|_{p\M p})_+$ which converges to $L_p(x)$. Applying \cite[Proposition~4.23(2)]{hiai2021} once more, as in the proof of Theorem~\ref{thm:order-free-allocation}, gives $\widehat L_p(\widehat x)=\sup_N\widehat{S_N}=\widehat{L_p(x)}$.

(ii) and (iii) Let $\mathcal H$ be the Hilbert space on which $\M$ acts, so that $p\M p$ acts on $p\mathcal H$. Since the $q_n$ are pairwise orthogonal with $\sum_nq_n=p$, for Borel sets $S\subseteq[0,\infty)$, the formula $E(S):=\sum_{n\in S\cap\N}q_n$ defines a projection-valued measure on $[0,\infty)$ with values in $p\M p$ and $E([0,\infty))=p$. Let $\Phi_p:=\int_{[0,\infty)}t~dE(t)$ be the corresponding positive self-adjoint operator on $p\mathcal H$. All its spectral projections $E(S)$ belong to $p\M p$, so $\Phi_p$ is affiliated with $p\M p$, and $\one_{\{n\}}(\Phi_p)=E(\{n\})=q_n$. Extend $h$ to a Borel function $h\colon[0,\infty)\to[0,\infty)$; since $E$ is concentrated on $\N$, the operator $H:=h(\Phi_p)=\int h~dE$ does not depend on the extension, and it is $\tau$-measurable, since $p\M p$ is finite. For $N\geq1$, let $P_N:=\sum_{n\leq N}|B_n|q_n\in(p\M p)_+$. Then $H-P_N=(h\one_{(N,\infty)})(\Phi_p)\geq0$, and \cite[(4.6), p.~50]{hiai2021} gives $\tau(H)=\sum_{n\geq1}|B_n|~\tau(q_n)$ and $\tau(H-P_N)=\sum_{n>N}|B_n|~\tau(q_n)$.
On the other hand, every branch $(n,\gamma)$ contributes $v_{n,\gamma}^*v_{n,\gamma}=q_n$ to $L_p(1)$, so $P_N$ is the partial sum of the series $L_p(1)$ over the branches $(n,\gamma)$ with $n\leq N$. By Proposition~\ref{prop:universal-kac-towers3}, $P_N\to L_p(1)$ in $L^1(p\M p,\tau|_{p\M p})$, and $\tau(L_p(1))=\tau(1)=1$. Hence $\sum_n|B_n|~\tau(q_n)=\lim_N\tau(P_N)=1$. Therefore $H\in L^1(p\M p,\tau|_{p\M p})_+$ and $\|H-P_N\|_1=\tau(H-P_N)\to0$. By uniqueness of $L^1$-limits, $R_p=L_p(1)=H=h(\Phi_p)$ and $\tau(R_p)=\sum_n|B_n|~\tau(q_n)=1$. Since $p\neq0$ and $\tau$ is faithful, $\tau(p)>0$, and dividing by $\tau(p)$ gives $\tau_p(R_p)=1/\tau(p)$.
\end{proof}

\begin{corollary}
\label{cor:universal-towers-ergodic}
Suppose that $\alpha$ is ergodic and preserves a faithful normal state $\varphi$ on the finite von Neumann algebra $\M$. Then $\varphi$ is tracial, every non-zero projection $p\in\M$ is a sweep-out projection, and Proposition~\ref{prop:universal-kac-towers1} and Theorem~\ref{thm:universal-kac-towers} apply with $\tau=\varphi$.
\end{corollary}

\begin{proof}
By Proposition~\ref{prop:invariant-state-is-trace}, $\varphi$ is a faithful normal tracial state, and the argument at the beginning of Section~\ref{sec:universal-towers} shows that every non-zero projection is a sweep-out projection.
\end{proof}

For $k\geq1$, let $H_k:=\one_{[k,\infty)}(R_p)$.

\begin{corollary}
\label{cor:universal-tower-tail}
For every $k\geq1$, the projection
$H_k$ is the sum of those
$q_n$ for which $|B_n|\geq k$, and
$\tau(H_k)\leq1/k$.
\end{corollary}

\begin{proof}
By Theorem~\ref{thm:universal-kac-towers}\textup{(iii)}, $R_p=h(\Phi_p)$, and the functional calculus gives
\[
H_k=(\one_{[k,\infty)}\circ h)(\Phi_p)=E(\{n\in\N:|B_n|\geq k\})=\sum_{n:~|B_n|\geq k}q_n,
\]
where $E$ is the spectral measure of $\Phi_p$ from the proof of Theorem~\ref{thm:universal-kac-towers}. For $t\geq0$, we have $k\one_{[k,\infty)}(t)\leq t$, so $kH_k\leq R_p$, and \cite[Corollary~4.21(1)]{hiai2021} gives $k\tau(H_k)\leq\tau(R_p)=1$.
\end{proof}

\begin{remark}
\label{rem:towers-allocation}
A partial-isometric Kac tower defines a partial-isometric $p$-allocation. For $\gamma\in\Gamma$, let $N_\gamma:=\{n\geq1:\gamma\in B_n\}$, and let $\rho_\gamma:=\sum_{n\in N_\gamma}\rho_{n,\gamma}$, $p_\gamma:=\sum_{n\in N_\gamma}q_n$ and $u_\gamma:=\sum_{n\in N_\gamma}v_{n,\gamma}$.
For distinct $n,n'\in N_\gamma$, the initial projections $q_n,q_{n'}$ of $v_{n,\gamma},v_{n',\gamma}$ are orthogonal, and so are their final projections $\alpha_\gamma(\rho_{n,\gamma}),\alpha_\gamma(\rho_{n',\gamma})$. Hence, for every vector $\xi$ and every finite set $K\subseteq N_\gamma$, we have $\|\sum_{n\in K}v_{n,\gamma}\xi\|^2=\sum_{n\in K}\|q_n\xi\|^2$. Therefore the series defining $u_\gamma$ converges strongly, and for every $\xi$, we have $\|u_\gamma\xi\|^2=\|p_\gamma\xi\|^2$, so that $u_\gamma^*u_\gamma=p_\gamma\leq p$, and the range of $u_\gamma$ is the closed linear span of the ranges of the projections $\alpha_\gamma(\rho_{n,\gamma})$, $n\in N_\gamma$, so that $u_\gamma u_\gamma^*=\alpha_\gamma(\rho_\gamma)$ by normality of $\alpha_\gamma$. By Definition~\ref{def:kac-tower}\textup{(ii)}, the $\rho_\gamma$ are pairwise orthogonal with $\sum_\gamma\rho_\gamma=1$. Thus $(\rho_\gamma)_{\gamma\in\Gamma}$ is a partial-isometric $p$-allocation with transport witness $(p_\gamma,u_\gamma)_{\gamma\in\Gamma}$, and Theorem~\ref{thm:order-free-allocation} applies to it. Moreover, for $n\in N_\gamma$, we have $u_\gamma q_n=v_{n,\gamma}$, and for $n\notin N_\gamma$, we have $u_\gamma q_n=0$. Consequently, $L_{\rho,p,u}(1)=\sum_{\gamma\in\Gamma}p_\gamma=\sum_{n\geq1}|B_n|~q_n=R_p$, where we rearranged a series of positive terms whose $L^1$-norms have finite sum. For $x\in\M$ and $n\geq1$, the compression $y\mapsto q_nyq_n$ is continuous on $L^1(p\M p,\tau|_{p\M p})$ by \cite[Lemma~4.35(2)]{hiai2021}, and it gives $q_nL_{\rho,p,u}(x)q_n=\sum_{\gamma\in B_n}v_{n,\gamma}^*\alpha_\gamma(x)v_{n,\gamma}$. Hence $L_p(x)=\sum_{n\geq1}q_nL_{\rho,p,u}(x)q_n$ is the block-diagonal part of $L_{\rho,p,u}(x)$ with respect to the decomposition $p=\sum_nq_n$, and the two operators coincide for $x=1$.
\end{remark}

In the following proposition, we show that, for abelian $\M$, Proposition~\ref{prop:universal-kac-towers1} and Theorem~\ref{thm:universal-kac-towers} recover the universal-shape theorem of Meyerovitch and Weiss \cite[Theorem~1.5]{mw24}.

Let $T\colon\Gamma\curvearrowright(X,\mathcal B,\mu)$ be a measure-preserving left action on a probability space, $\M=L^\infty(X,\mu)$, $\tau(f)=\int_Xf~d\mu$ and $\alpha_\gamma(f)=f\circ T_\gamma^{-1}$. Let $A\in\mathcal B$ be a sweep-out set, that is, the orbit of almost every point of $X$ meets $A$, so that $\bigcup_{\gamma\in\Gamma}T_\gamma(A)=X$ modulo $\mu$, and let $p=\one_A$. For a measurable partition $A=\biguplus_{n\geq1}A_n$, let $\phi_A\colon A\to\N$ be given by $\phi_A(z)=n$ for $z\in A_n$.

\begin{proposition}
\label{prop:commutative-universal-towers}
Let $(B_n)_{n\geq1}$ be a sequence of finite subsets of $\Gamma$ containing $e$.
\begin{enumerate}
\item For every partial-isometric Kac tower for $(\alpha,p)$ with shapes $(B_n)_{n\geq1}$, there is a measurable partition $A=\biguplus_{n\geq1}A_n$ such that $q_n=\one_{A_n}$, $\rho_{n,\gamma}=\one_{T_\gamma^{-1}(A_n)}$ for every branch $(n,\gamma)$, and
\begin{equation}
\label{eq:MW-universal-cover}
X=\biguplus_{n\geq1}\biguplus_{\gamma\in B_n}T_\gamma^{-1}(A_n)~\text{modulo }\mu.
\end{equation}
Conversely, for every measurable partition $A=\biguplus_{n\geq1}A_n$ satisfying \eqref{eq:MW-universal-cover}, the elements $q_n=\one_{A_n}$, $\rho_{n,\gamma}=\one_{T_\gamma^{-1}(A_n)}$ and $v_{n,\gamma}=q_n$ form a partial-isometric Kac tower.
\item For a partial-isometric Kac tower and a partition as in \textup{(i)}, we have $R_p=|B_{\phi_A}|$ almost everywhere on $A$, and $\int_A|B_{\phi_A(z)}|~d\mu(z)=1$. For every $f\in L^1(X,\mu)$, we have $\int_AL_p(f)~d\mu=\int_Xf~d\mu$ and, for almost every $z\in A$, we have $L_p(f)(z)=\sum_{\gamma\in B_{\phi_A(z)}}f(T_\gamma^{-1}(z))$.
\end{enumerate}
In particular, if $(B_n)_{n\geq1}$ is the sequence of Proposition~\ref{prop:universal-kac-towers1}, then for every probability-preserving action of $\Gamma$ and every sweep-out set $A$, there is a measurable map $\phi_A\colon A\to\N$ such that \eqref{eq:MW-universal-cover} holds with $A_n=\phi_A^{-1}(n)$ and $\int_A|B_{\phi_A(z)}|~d\mu(z)=1$. This is \cite[Theorem~1.5]{mw24}.
\end{proposition}

\begin{proof}
We use that, for $E\in\mathcal B$, we have $\alpha_\gamma(\one_E)=\one_{T_\gamma(E)}$, and that equivalent projections in the abelian algebra $\M$ are equal (Remark~\ref{rem:direct-allocation}).

(i) Let $(q_n,\rho_{n,\gamma},v_{n,\gamma})$ be a partial-isometric Kac tower, and write $q_n=\one_{A_n}$ and $\rho_{n,\gamma}=\one_{E_{n,\gamma}}$. By Definition~\ref{def:kac-tower}\textup{(iii)}, $\alpha_\gamma(\rho_{n,\gamma})=v_{n,\gamma}v_{n,\gamma}^*=v_{n,\gamma}^*v_{n,\gamma}=q_n$, that is, $T_\gamma(E_{n,\gamma})=A_n$, or equivalently $E_{n,\gamma}=T_\gamma^{-1}(A_n)$, modulo null sets. By Definition~\ref{def:kac-tower}\textup{(i)}, the sets $A_n$ partition $A$ modulo null sets; modifying them on null sets, we obtain a partition of $A$. By Definition~\ref{def:kac-tower}\textup{(ii)}, the sets $T_\gamma^{-1}(A_n)$ satisfy \eqref{eq:MW-universal-cover}. Conversely, let $A=\biguplus_nA_n$ be a measurable partition satisfying \eqref{eq:MW-universal-cover}, and let $q_n=\one_{A_n}$, $\rho_{n,\gamma}=\one_{T_\gamma^{-1}(A_n)}$ and $v_{n,\gamma}=q_n$. Properties (i) and (ii) of Definition~\ref{def:kac-tower} follow from the two partitions, $v_{n,\gamma}^*v_{n,\gamma}=v_{n,\gamma}v_{n,\gamma}^*=q_n=\alpha_\gamma(\rho_{n,\gamma})$ gives (iii), and $e\in B_n$ gives $\rho_{n,e}=q_n=v_{n,e}$, which is (iv).

(ii) For $f\in L^\infty(X,\mu)$, since $\M$ is abelian, we have
\[
v_{n,\gamma}^*\alpha_\gamma(f)v_{n,\gamma}=v_{n,\gamma}^*v_{n,\gamma}\alpha_\gamma(f)=\one_{A_n}\cdot(f\circ T_\gamma^{-1}).
\]
Both sides are continuous in $f$ for the $L^1$-norm, by \cite[Lemma~4.35(2)]{hiai2021} and Proposition~\ref{prop:universal-kac-towers2}, so this identity also holds for $f\in L^1(X,\mu)$. For $z\in A_n$, the only branches contributing at $z$ are $(n,\gamma)$ with $\gamma\in B_n$, and their contributions are $f(T_\gamma^{-1}(z))$. Since $B_n$ is finite, the partial sums of $L_p(f)$ along an enumeration of the branches converge almost everywhere on $A$ to $\sum_{\gamma\in B_{\phi_A(z)}}f(T_\gamma^{-1}(z))$. By Proposition~\ref{prop:universal-kac-towers3}, they also converge in $L^1(X,\mu)$ to $L_p(f)$, so a subsequence converges almost everywhere to $L_p(f)$, which proves the formula for $L_p(f)$. Taking $f=1$ gives $R_p=L_p(1)=|B_{\phi_A}|$ almost everywhere on $A$. The two integral identities are $\tau(L_p(f))=\tau(f)$ (Proposition~\ref{prop:universal-kac-towers3}) and $\tau(R_p)=1$ (Theorem~\ref{thm:universal-kac-towers}).

The last assertion follows from Proposition~\ref{prop:universal-kac-towers1}, applied to $(L^\infty(X,\mu),\tau,\alpha,\one_A)$, together with (i) and (ii); here $\one_A$ is a sweep-out projection, because $\bigvee_{\gamma\in\Gamma}\alpha_\gamma(\one_A)=\one_{\bigcup_\gamma T_\gamma(A)}=1$.
\end{proof}

\begin{remark}
\label{rem:tower-choice}
The sequence $(B_n)_{n\geq1}$ of Proposition~\ref{prop:universal-kac-towers1} depends only on $\Gamma$, but a partial-isometric Kac tower with these shapes is not unique, and $R_p=L_p(1)$ depends on it. For instance, in the example of Remark~\ref{rem:allocation-identity-scope}, let $n_1,n_2,n_3$ be indices with $B_{n_1}=\{0,1\}$, $B_{n_2}=\{0\}$ and $B_{n_3}=\{0,2\}$. The partitions of $A=\{0,1\}$ given by $A_{n_1}=\{0\}$, $A_{n_2}=\{1\}$, and by $A_{n_2}=\{0\}$, $A_{n_3}=\{1\}$, all other $A_n$ being empty, satisfy \eqref{eq:MW-universal-cover}, since $T_1^{-1}(\{0\})=\{2\}$ and $T_2^{-1}(\{1\})=\{2\}$. By Proposition~\ref{prop:commutative-universal-towers}, they define partial-isometric Kac towers with $R_p=2\cdot\one_{\{0\}}+\one_{\{1\}}$ and $R_p=\one_{\{0\}}+2\cdot\one_{\{1\}}$, respectively; in both cases $\tau_p(R_p)=1/\tau(p)=3/2$, in accordance with Theorem~\ref{thm:universal-kac-towers}. In contrast, for $\Gamma=\Z$ the return time $T_p$ of Theorem~\ref{thm:quantum-return-time} is determined by the projections $\alpha^k(p)$, $k\geq1$. In this example, $\alpha=\alpha_1$ is given by $\alpha(f)=f\circ S$ with $S(y)=y-1$, and Proposition~\ref{prop:commutative-Z} gives $T_p=r_A=2\cdot\one_{\{0\}}+\one_{\{1\}}$. Thus the first tower satisfies $R_p=T_p$, whereas the second does not.
\end{remark}

\section*{Declaration on the use of AI tools}
The authors proved all the results for $\Z$-actions in Section~\ref{sec:Z-action} in January~2026, without AI assistance. The allocation and fibre-sum motivation, and the aim of extending it to actions on von Neumann algebras, are due to the authors. During the development of this extension, ChatGPT-5.6 sol contributed the mathematical construction used to overcome the obstruction described at the beginning of Section~\ref{subsec:coupling}: one first passes to the center-valued dimensions $\Theta(p)$ and $\Theta(p_j)$, constructs central positive elements $c_F$ with $\sum_Fc_F=\Theta(p)$ and $\sum_{F\ni g_j}c_F=\Theta(p_j)$ (Lemma~\ref{lem:finite-shape-threshold-coupling}), and then realizes these dimensions by projections and partial isometries (Proposition~\ref{prop:universal-kac-towers1}). This contribution is mathematical and is not merely language editing.

\bibliographystyle{amsplain}
\bibliography{newreferences}

\end{document}